\documentclass[aos,seceqn,citesort,dvips]{arximspdf}
\usepackage{mathbh}
\usepackage{graphicx}

\doi{10.1214/10-AOS806}
\volume{38}
\issue{5}
\pubyear{2010}
\firstpage{2781}
\lastpage{2822}

\makeatletter

\newtheorem{Th}{Theorem}
\newtheorem{Prop}{Proposition}
\newtheorem{Lemma}{Lemma}

\newproclaim{remark}{Remark}

\makeatother

\begin{document}
\begin{frontmatter}

\title{Adaptive estimation for Hawkes processes; application to genome
analysis\thanksref{T1}}
\runtitle{Adaptive estimation for Hawkes' processes}

\thankstext{T1}{Supported by the French Agence Nationale
de la Recherche (ANR), under Grant ATLAS
(JCJC06\_137446) ``From Applications to Theory in Learning and
Adaptive Statistics.''}

\begin{aug}
\author[A]{\fnms{Patricia} \snm{Reynaud-Bouret}\corref{}\ead[label=e1]{reynaudb@unice.fr}} and
\author[B]{\fnms{Sophie} \snm{Schbath}\thanksref{t2}\ead[label=e2]{Sophie.Schbath@jouy.inra.fr}}
\runauthor{P. Reynaud-Bouret and S. Schbath}
\affiliation{CNRS, Universit\'e de Nice Sophia-Antipolis and Institut
National de~la~Recherche~Agronomique}
\address[A]{Laboratoire J. A. Dieudonn\'e\\
U.M.R. C.N.R.S. 6621 \\
Universit\'e de Nice Sophia-Antipolis\\
Parc Valrose\\
06108 Nice Cedex 2\\
France\\
\printead{e1}}
\address[B]{Institut National de la Recherche Agronomique \\
Unit\'e Math\'ematique, Informatique et G\'enome \\
Domaine de Vilvert \\
F-78352 Jouy-en-Josas Cedex \\
France\\
\printead{e2}}
\end{aug}

\thankstext{t2}{Most of this work has been done while
P. Reynaud-Bouret was working at the Ecole Normale Sup\'erieure de
Paris (DMA, UMR 8553).}

\received{\smonth{3} \syear{2009}}
\revised{\smonth{12} \syear{2009}}

%
\begin{abstract}
The aim of this paper is to provide a new method for the detection of
either favored or avoided distances between genomic events along DNA
sequences. These events are modeled by a Hawkes process. The biological
problem is actually complex enough to need a nonasymptotic penalized
model selection approach. We provide a theoretical penalty that
satisfies an oracle inequality even for quite complex families of
models. The consecutive theoretical estimator is shown to be adaptive
minimax for H\"olderian functions with regularity in $(1/2,1]$: those
aspects have not yet been studied for the Hawkes' process. Moreover,
we introduce an efficient strategy, named \textit{Islands}, which is not
classically used in model selection, but that happens to be
particularly relevant to the biological question we want to answer.
Since a multiplicative constant in the theoretical penalty is not
computable in practice, we provide extensive simulations to find a
data-driven calibration of this constant. The results obtained on real
genomic data are coherent with biological knowledge and eventually
refine them.
\end{abstract}

%
\begin{keyword}[class=AMS]
\kwd[Primary ]{62G05}
\kwd{62G20}
\kwd[; secondary ]{46N60}
\kwd{65C60}.
\end{keyword}
\begin{keyword}
\kwd{Hawkes process}
\kwd{model selection}
\kwd{oracle inequalities}
\kwd{data-driven penalty}
\kwd{minimax risk}
\kwd{adaptive estimation}
\kwd{unknown support}
\kwd{genome analysis}.
\end{keyword}

\end{frontmatter}

\section{Introduction}\label{sec1}

Modeling the arrival times of a particular event on the real
line is a common problem in time series theory. In this paper, we
deal with a very similar
but rarely addressed problem: modeling
the process of the occurrences of a particular event along a
discrete sequence, namely a DNA sequence. Such events could be, for
instance, any given DNA patterns, any genes or any other biological signals
occurring along genomes. A huge literature exists on the statistical
properties of pattern occurrences along random sequences \cite{ReiS00} but
our current approach is different. It consists in directly modeling the
point process of the occurrences of any kind of events and it is not
restricted to pattern
occurrences.
Our aim is to characterize the dependence, if any,
between the event occurrences by pointing out either favored or
avoided distances
between them, those distances being significantly larger than the
classical memory used in the quite popular Markov chain model for
instance. At this scale, it is more interesting to use a continuous
framework and see occurrences as points. A very interesting model for
this purpose is the
Hawkes process \cite{fado}.

In the most basic self-exciting model, the Hawkes process $(N_t)_{t\in
\mathbb{R}}$ is defined by its intensity, which satisfies
%
%
\begin{equation}
\label{class}
\lambda(t)= \nu+\int_{-\infty}^{t^-} h(t-u)\,dN_u,
\end{equation}
where $\nu$ is a positive parameter, $h$ a nonnegative function with
support on $\mathbb{R}^+$ and $\int h<1$ and where $dN_u$ is the point measure
associated to the process. The interested reader shall find in Daley
and Vere-Jones' book \cite{daleyVerejones} the main definitions,
constructions and models related to point processes in general and
Hawkes processes in particular [see, e.g., Examples 6.3(c) and
7.2(b) therein].

The intensity $\lambda(t)$ represents the probability to have an
occurrence at position $t$ given all the past. In this sense,
(\ref{class}) basically means that there is a constant rate $\nu$ to
have a spontaneous occurrence at $t$ but that also all the previous
occurrences influence the apparition of an occurrence at $t$. For
instance, an occurrence at $u$ increases the intensity by $h(t-u)$. If
the distance $d=t-u$ is favored, it means that $h(d)$ is really large:
having an occurrence at $u$ significantly increases the chance of
having an occurrence at $t$. The intensity given by (\ref{class}) is
the most basic case, but variations of it
enable us to model self-inhibition, which happens when one allows $h$
to take negative values (see Section \ref{sec24}) and, in the
most general case, to model interaction with another type of event. The
drawback is that, by definition, the Hawkes process is defined on an
ordered real line (there is a past, a present and a future). But a
strand of DNA itself has a direction, a fact that makes our approach
quite sensible.

The Hawkes model has been widely used to model
the
occurrences of earthquake \cite{VeO82}. In this set-up and even for
more general counting processes, the statistical inference usually
deals with maximum likelihood estimation \cite{OgA82,Oz79}.
This approach has been applied to genome analysis: in
a previous work \cite{fado}, Gusto and Schbath's method, named FADO,
uses maximum likelihood estimates of the coefficients of $h$ on a
Spline basis coupled with an AIC criterion to select the set of
equally spaced knots.

On one hand, the FADO procedure is quite effective---it can manage
interactions between two types of events and self excitation or
inhibition, that is, it works in the most general Hawkes process framework
and produces smooth estimates. However, there are several
drawbacks. From a theoretical point of view, AIC criterion is proved
to select the right set of knots if first, there exists a true set of
knots, and then if the family of possible knots is held fixed whereas
the length of the observed sequence of DNA tends to infinity. Moreover,
from a practical point of view, the criterion seems to behave very
poorly when a lot of possible sets of knots with the same cardinality
are in competition \cite{gaelle}. FADO has been implemented with
equally spaced knots for this reason. Finally, it heavily depends on an
extra knowledge of the support of the function $h$. In practice, we
have to input the maximal size of the support,
say 10,000 bases, in the
FADO procedure. Consequently the FADO estimate is a spline
function based on knots that are equally spaced on
$[0,10\mbox{,}000]$. If this maximal size is too large,
the estimate of $h$ will probably be small with
some fluctuations but not null until the end of the interval, whereas
it should be null before (see Figure \ref{avecfado} in Section \ref{sec5}).

On the other hand, our feeling is that if interaction exists, say
around the distance $d=500$ bases, the function $h$ to estimate should
be really large, around $d=500$, and if there is no biological reason
for any other interaction, then $h$ should be null anywhere else.

One way to solve this problem of estimation is to use model selection
but in its nonasymptotic version. Ideally, if the work of Birg\'e and
Massart in \cite{gauss} was not restricted to the Gaussian case but if
it also provides results for the Hawkes model then it should enable
us to find a way of selecting an irregular set of knots with
complexity that may grow if the length of the observed sequence
becomes larger. The question of the knowledge of the support never
appears in Birg\'e and Massart's work because there is not such a
question in a Gaussian model, but one could imagine that their way of
selecting sparse models should enable us to select a sparse support
too.

However, we are not in an ideal world where a white noise model and
Hawkes model are equivalent (even heuristically), so there is no way
to guess the right way of penalizing in our situation. So the
purpose of this article is to provide a first attempt at constructing a
penalized model
selection in a nonasymptotic way for the Hawkes model. This paper
consists in both practical methods for estimating $h$ that lie on
theoretical evidences and also in new theoretical results such as
oracle inequalities or adaptivity in the minimax sense. Note that,
to our knowledge, the minimax aspects of the Hawkes model have not yet
been considered.

Accordingly, we restrict ourselves to a simpler case than the FADO
procedure. First, we focus on the self-exciting model [i.e., the one
given by (\ref{class}), where $h$ is assumed to be nonnegative], but we
would at least like that the final
estimator remains computable in case of self-inhibition. Then we do
not use maximum likelihood estimators since they are not easily
handled by model selection procedures, at least from a theoretical
point of view. So we provide in this paper theoretical
results for penalized projection estimators (i.e., least square
estimators) and not for penalized maximum likelihood estimators (see
Chapter 7 of \cite{stflourpascal} for a complete comparison of both
contrasts in the density setting from a model selection point of view).
Finally, for technical reasons, we only
deal with piecewise constant estimators.
Once all those restrictions are done, the gap between the theoretical
procedure and the practical procedure is consequently reduced to a
practical calibration problem of the multiplicative constants.

Since the Hawkes processes are quite popular for modeling
earthquakes, financial, or economical data, we try to keep a general
formalism in most of the
sequel (except in the biological applications part). Consequently, our
method could be applied to many other type
of data.

In Section \ref{sec2}, we define the notation and the different
families of
models. Section \ref{sec3} states first a nonasymptotic result for the
projection estimators, since up to our knowledge, these estimators were
not yet studied. Then Section \ref{sec3} gives a theoretical penalty that
enables us to select
a good estimator in a family of projection estimators. Indeed, we prove
that our penalized projection estimator satisfies an oracle
inequality, hence proving by that result that our estimator is as good
as the
best projection estimator in the family up to some multiplicative term. However,
the multiplicative constant in the theoretical penalty is not
computable in practice. As a consequence, Section \ref{sec4}
provides simulations which validate a calibration method that seems
to work well from a practical point of view. Then in Section~\ref{sec5} we
apply this method to DNA data. The results match biological evidences
and refine them. Section \ref{sec6} details the adaptive and minimax
properties of our estimators. Section \ref{sec7} is dedicated to more technical
results that are at the origin of the ones stated in Section \ref{sec3}.
Sketch of proofs can be found in Section \ref{sec8}: the interested
reader shall
find details of those proofs in \cite{sopacomplet}.

\section{Framework}\label{sec2}
Let $(N_t)_{t}$ be a stationnary Hawkes process on the real line
satisfying (\ref{class}). We assume that $h$ has a bounded support
included in $(0,A]$ where $A$ is a known positive real number and that
%
%
\begin{equation}\label{defp}
p:=\int_0^A h(u) \,du
\end{equation}
satisfies $p<1$. This condition guarantees the existence of a stationary
version of the process (see \cite{HaO74}).
Let us remark that, for the DNA applications we have in mind,
$A$ is quite known because it corresponds to a maximal distance
from which it is no longer reasonable to consider a linear
interaction between two genomic locations. If there may exist some
interaction at longer distances, then it should certainly imply the 3D
structure of DNA.

We observe the stationary Hawkes process $(N_t)_{t}$ on an interval
$[-A,T]$, where $T$ is a positive real number. Typically $T$ should
be significantly larger than $A$. Using this observation, we want to
estimate
%
%
\begin{equation}
\label{defs}
s=(\nu,h),
\end{equation}
assumed to be in
%
%
\begin{eqnarray}\label{l2}
\mathbb{L}^{2} &=& \biggl\{f=(\mu,g)\dvtx g \mbox{ with support in }
(0,A],\nonumber\\[-8pt]\\[-8pt]
&&\hspace*{15pt}
\Vert f \Vert^2=\mu^2+\int_0^Ag^2(x)\,dx <+\infty\biggr\}.\nonumber
\end{eqnarray}
The introduction of this Hilbert space is related to the fact that we
want to use least square estimators.

With these constraints on $h$, we can note that (\ref{class}) is
equivalent to
%
%
\begin{equation}
\label{intensiteclass}
\lambda(t)=\nu+ \int_{t-A}^{t^-} h(t-u) \,dN_u.
\end{equation}
Now, we can introduce intensity candidates:
for all $f=(\mu,g)$ in $\mathbb{L}^{2}$, we define
%
%
\begin{equation}
\label{psi}
\Psi_{f}(t):= \mu+ \int_{t-A}^{t^-} g(t-u) \,dN_u.
\end{equation}
In particular, note that $\Psi_{s}(t)=\lambda(t)$. A good intensity candidate
should be a $\Psi_{f}(\cdot)$ that is close to $\Psi_{s}(\cdot)$. The least-square
contrast is consequently defined for all $f$ in $\mathbb{L}^{2}$ by
%
%
\begin{equation}
\label{contraste}
\gamma_T(f):= -\frac{2}{T} \int_{0}^{T} \Psi_{f}(t) \,dN_t + \frac{1}{T}
\int_{0}^{T}
\Psi_{f}(t)^2 \,dt.
\end{equation}
As we will see in Lemma \ref{contra}, this really defines a contrast,
in the statistical sense.
Indeed, taking the compensator of the previous formula leads to
\[
-\frac{2}{T} \int_{0}^{T} \Psi_{f}(t) \Psi_{s}(t) \,dt + \frac{1}{T}
\int_{0}^{T}
\Psi_{f}(t)^2 \,dt.
\]
Let us consider the last integral in the previous equation:
%
%
\begin{equation}
\label{DT}
D^2_T(f):= \frac{1}{T} \int_{0}^{T} \Psi_{f}(t)^2 \,dt.
\end{equation}
Lemma \ref{norme} proves that $
D^2_T(\cdot)$
defines a quadratic form on $\mathbb{L}^{2}$ such
that
%
%
\begin{equation}
\label{normD}
\Vert f \Vert_D:=\sqrt{\mathbb{E}(D^2_T(f))}
\end{equation}
is a quadratic norm on $\mathbb{L}^{2}$, equivalent to
$\Vert f \Vert$ [see (\ref{l2})].
In this sense, we can see $\gamma_T(f)$ as an empirical version of
$\Vert f-s \Vert_D^2-\Vert s \Vert_D^2$, which is quite classical for a
least-square contrast (see the density set-up, e.g., in
\cite{stflourpascal}).

\subsection{Projection estimator}\label{sec21}
Let $m$ be a set of disjoint intervals of $(0,A]$. In the sequel, $m$ is
called a model and $|m|$ denotes the number of intervals in
$m$. One can think of $m$ as a partition of $(0,A]$ but there are
other interesting cases as we will see later.
Let $S_m$ be the vectorial space of $\mathbb{L}^{2}$ defined by
%
%
\begin{equation}
\label{modele}
\quad S_m= \biggl\{ f=(\mu,g) \in\mathbb{L}^{2} \mbox{ such that } g =
\sum_{I \in m} a_I\frac{\mathbh{1}_{I}}{\sqrt{\ell(I)}} \mbox{
with }
(a_I)_{I\in m}\in\mathbb{R}^m \biggr\},\hspace*{-28pt}
\end{equation}
where $\ell(I)=\int\mathbh{1}_I \,dt$. We say that $g$ in the above equation
is constructed on the model $m$. Conversely, if $g$ is a piecewise
constant function, remark that we can define a resulting model $m$ by
the set of intervals where $g$ is constant but nonzero and a resulting
partition by the set of intervals where $g$ is constant.
The projection estimator, $\hat{s}_m$, is the least square estimator of
$s$ defined by
%
%
\begin{equation}
\label{chaps}
\hat{s}_m:= \mathop{\arg\min}_{f \in S_m} \gamma_T(f).
\end{equation}
Of course the estimator $\hat{s}_m$ heavily depends on the choice of
the model $m$. That is the main reason for trying to select it in a
data driven way.
Model selection intuition usually relies on a bias-variance
decomposition of the risk of $\hat{s}_m$.
So let us define $s_m$ as the orthogonal projection for
\mbox{$\Vert\cdot \Vert$} of $s$ on $S_m$. Then $\hat{s}_m$ is
a ``good'' estimate of $s_m$, since $\gamma_T(f)$ is an approximation of
$\Vert f-s \Vert_D^2-\Vert s \Vert_D$. We cannot prove that it is an unbiased
estimate, but the intuition applies. So the bias can be more or less
identified as $\Vert s-s_m \Vert^2$. This is the approximation error
of the
model $m$ with respect to $s$.
As we will see in Proposition \ref{surunmodele} and the consecutive
comments, one can actually prove that
\[
\mathbb{E}(\Vert s-\hat{s}_m \Vert^2)\simeq C_T \biggl[\Vert s-s_m \Vert
^2+\frac
{|m|}{T} \biggr],
\]
where $C_T$ is a positive quantity that slowly varies with $T$. So the
variance or stochastic error may be identified as $|m|/T$. We recover a
bias-variance decomposition where the bias decreases and the variance
increases. Finding a model $m$ in a data driven way that almost
minimizes the previous equation is the main goal of model selection.
However, there is no precise shape for the quantity $C_T$.
We consequently use the most general form of penalization in the sequel.

\subsection{Penalized projection estimator}\label{sec22}
Let $\mathcal{M}_{T}$ be a family of sets of disjoint intervals of $(0,A]$
(i.e., a family of possible models). We denote by $\#\{\mathcal
{M}_{T}\}$ the total
number of models. We define the penalty (or penalty function) by
$\operatorname{pen}
\dvtx\mathcal{M}_{T}\to\mathbb{R}^+$ and we select a model by
minimizing the following criterion:
%
%
\begin{equation}
\label{critpen}
\hat{m}:=\mathop{\arg\min}_{m\in\mathcal{M}_{T}} [\gamma_T(\hat{s}_m) +
\operatorname{pen}(m)].
\end{equation}
Then the penalized projection estimator is defined by
%
%
\begin{equation}
\label{estpen}
\tilde{s}=
(\tilde{\nu},\tilde{h})=\hat{s}_{\hat{m}}.
\end{equation}
The main problem is now to find a function $\operatorname{pen}\dvtx
\mathcal{M}_{T}\to\mathbb{R}^+$ that
guarantees that
%
%
\begin{equation}\label{oracleessai}
\Vert s-\tilde{s} \Vert^2\leq {C \inf_{m\in\mathcal{M}_{T}}} \Vert
s-\hat{s}_m \Vert^2
\end{equation}
and this either with high probability or in expectation, up to some
small residual term and up to some multiplicative term $C$ that could
slightly increase with $T$. The previous equation (\ref{oracleessai})
is an oracle inequality. If
this oracle inequality holds, this will mean that we can select a model
$\hat{m}$, and consequently a projection estimator $\tilde{s}=\hat
{s}_{\hat{m}}$, that is almost as good as the best estimator in the
family of the $\hat{s}_m$'s---whereas this best estimator cannot be
guessed without knowing $s$. Of course this would tell us nothing if
the projection estimators themselves, that is, the $\hat{s}_m$'s, are
not sensible. The next section precisely states the properties of the
projection estimator and the oracle inequality satisfied by the
penalized projection estimator.
To conclude Section \ref{sec2}, we precise the different families of
models we
would like to use and we precisely explain what self-inhibition means
in our model.

\subsection{Strategies}\label{sec23}
A strategy refers to the choice of the family of models $\mathcal{M}_{T}$.
In the sequel, a partition $\Gamma$ of $(0,A]$ should be understood as
a set of disjoint intervals of $(0,A]$ such that their union is the
whole interval $(0,A]$. A regular partition is such that all its
intervals have the same length.
We say that a model $m$ is written on $\Gamma$ if all the extremities
of the intervals in $m$ are also extremities of intervals in $\Gamma$.
For instance if $\Gamma=\{(0,0.25],(0.25,0.5],(0.5 0.75],(0.75,1]\}$
then $\{(0,0.25],(0.25,1]\}$ or $\{(0,0.25],\break(0.75,1]\}$ are models
written on $\Gamma$.
Now let us give some examples of families $\mathcal{M}_{T}$.
Let $J$ and $N$ be two positive integers.

\textit{Nested strategy}. Take $\Gamma$ a dyadic regular partition (i.e.,
such that $|\Gamma|=2^J$). Then take $\mathcal{M}_{T}$ as the set of
all dyadic
regular partitions of $(0,A]$ that can be
written on $\Gamma$, including the
void set. In particular, note that $\#\{\mathcal{M}_{T}\}=J+2$. We say
that this
strategy is nested since for any pair of partitions in this family, one
of them is always written on the other one.

\textit{Regular strategy}. Another natural strategy is to look at all the
regular partitions of $(0,A]$ until some finest partition of cardinal
$N$. That is to say that one has exactly one model with cardinality
$k$ for each $k$ in $\{0,\ldots,N\}$.
Here $\#\{\mathcal{M}_{T}\}=N+1$.

\textit{Irregular strategy}.
Assume now that we know that $h$ is piecewise constant on $(0,A]$ but
that we
do not know where the cuts of the resulting partition are. We can
consider $\Gamma$ a regular partition such that $|\Gamma|=N$ and then
consider $\mathcal{M}_{T}$ the set of all possible partitions written
on $\Gamma$,
including the void set. In this case, $\#\{\mathcal{M}_{T}\}\simeq2^N$.

\textit{Islands strategy}. This last strategy has been especially designed
to answer our biological problem. We think that $h$ has a very localized
support. The interval $(0,A]$ is really large and in fact $h$ is
nonzero on a really
smaller interval or
a union of really smaller intervals: the resulting model is sparse. We
can consider $\Gamma$ a regular partition such that $|\Gamma|=N$ and
then consider $\mathcal{M}_{T}$ the set of all the
subsets of $\Gamma$. A typical $m$ corresponds to a vectorial space $S_m$
where the functions $g$ are zero on $(0,A]$ except on some disjoints intervals
which look like several ``islands.'' In this case, $\#\{\mathcal
{M}_{T}\}=2^N$.

%
%
\begin{figure}

\includegraphics{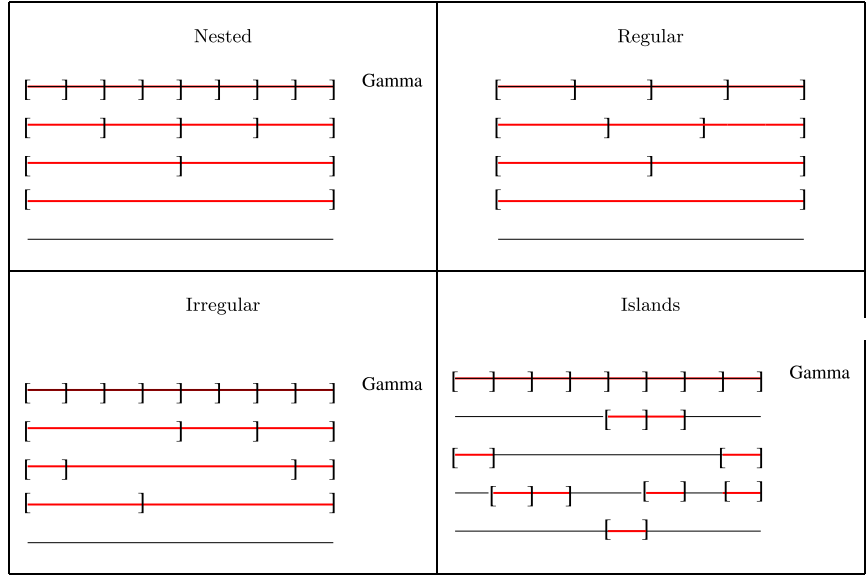}

\caption{On each line, one can find a model by looking at
the collection of red intervals between ``['' or ``].'' For the \textup{Nested
strategy}, here are all the models for $J=3$. For the \textup{Regular
strategy}, here are all the models for $N=4$. For the \textup{Irregular}
and \textup{Islands strategies}, these are just some examples of models
in the family with $N=8$.}\label{strategies}
\end{figure}

Figure \ref{strategies} gives some more visual examples of the
different strategies.

\subsection{Self-inhibition}\label{sec24}
The self-interaction can be modeled in a more general way by a process
whose intensity is given by
%
%
\begin{equation}
\label{sineg}
\lambda(t)= \biggl(\nu+\int_{-\infty}^{t^-} h(t-u)\,dN_u \biggr)_+,
\end{equation}
where $h$ may now be negative. We have taken the positive part to
ensure that the intensity remains positive. Then the condition $\int
|h|<1$ is
sufficient to ensure the existence of a stationary version of the process
(see \cite{BrM96}). When $h(d)$ is strictly positive there is a
self-excitation at distance $d$. When $h(d)$ is strictly negative,
then there is a self-inhibition. It is more or less the same
interpretation as above [see (\ref{class})] except that now all the
previous occurrences are voting whether they ``like'' or ``dislike''
to have a new occurrence at position $t$. If this process is not
studied in this paper from a theoretical point of view because of
major technical issues (except in the remarks following Theorem
\ref{oracleineqenproba}), note that however our projection estimators,
$\hat{s}_m$,
and penalized projection estimators, $\tilde{s}$, do not take the sign
of $g$
or $h$ into account for being computed. That is the reason why we will
use our
estimators, even in this case, for the numerical results.

Finally, we use in the sequel the notation $\lozenge$ which represents a
positive function of the parameters that are written in indices. Each
time $\lozenge_\theta$ is written in some equation, one
should understand that there exists a positive function of $\theta$
such that the equation holds. Therefore, the values of $\lozenge_\theta$
may change from line to line and even change in the same
equation. When no index appears, $\lozenge$ represents a positive
absolute constant.

\section{Main results}\label{sec3}
For technical reasons, we are not able to carefully control the
behavior of the projection estimators if $\nu$ tends to 0
or to infinity, but also if $p$ [see (\ref{defp})] tends to 1: in
such cases, the number of points in the process is either exploding or
vanishing. Consequently, the theoretical results are proved within a
subset of $\mathbb{L}^{2}$.
Let us define for all real numbers $H >0$, $\eta> \rho>0$, $1>P>0$, the
following
subset of $\mathbb{L}^{2}$:
\[
\mathcal{L}_{H,P}^{\eta,\rho}= \biggl\{f=(\mu,g)\in\mathbb{L}^{2} /
\mu\in[\rho,\eta] , g(\cdot) \in[0,H] \mbox{ and } \int_0^A g(u) \,du
\leq
P \biggr\}.
\]
If we know that $s$ belongs to $\mathcal{L}_{H,P}^{\eta,\rho}$ and if
we know the parameters $H,\eta$ and $\rho$, then it is reasonable to
consider the clipped projection estimator, $\bar{s}_m$. If we denote
the projection estimator $\hat{s}_m=(\hat{\nu}_m,\hat{h}_m)$, then
$\bar
{s}_m=(\bar{\nu}_m,\bar{h}_m)$ is given, for all positive $t$, by
%
%
\begin{equation}\label{clipunmodele}
\left\{\begin{array}{l}
\bar{\nu}_m = \cases{\hat{\nu}_m, &\quad if $\rho\leq\hat{\nu}_m
\leq
\eta$,\cr
\rho, &\quad if $\hat{\nu}_m <\rho$,\cr
\eta, &\quad if $\hat{\nu}_m >\eta$,}
\\[16pt]
\bar{h}_m (t) = \cases{\hat{h}_m(t), &\quad if $0\leq\hat{h}_m (t)\leq H$,
\cr
0, &\quad if $\hat{h}_m(t) <0$,\cr
H, &\quad if $\hat{h}_m(t) >H$.}
\end{array}\right.
\end{equation}
Note that $\bar{s}_m$, the clipped version of $\hat{s}_m$, is only
designed for theoretical purpose. Whereas $\hat{s}_m$ may be computed
even for possibly negative $h$, the computation of $\bar{s}_m$ does not
make sense in this more general framework.
For the clipped projection estimator, we can prove the following result.
\begin{Prop}
\label{surunmodele}
Let $(N_t)_{t\in\mathbb{R}}$ be a Hawkes process with intensity given
by $\Psi_{s}(\cdot)$. Let $m$
be a model written on $\Gamma$ where $\Gamma$
is a regular partition of $(0,A]$ such that
%
%
\begin{equation}\label{condgamma}
|\Gamma|\leq\frac{\sqrt{T}}{(\log T)^3}.
\end{equation}
Then if
$s$ belongs to $\mathcal{L}_{H,P}^{\eta,\rho}$, the clipped projection
estimator on the model $m$
satisfies
\[
\mathbb{E}(\Vert\bar{s}_m-s \Vert^2)\leq\lozenge_{H,P,\eta,\rho,A}
\biggl[\Vert s_m-s \Vert^2+(|m|+1)\frac{\log
T}{T} \biggr].
\]
\end{Prop}

This result is a control of the risk of the clipped projection
estimator on
one model. A first interpretation is to assume that $s$ belongs to
$S_m$. In this case, if $m$ is fixed whereas $T$ tends to infinity,
Proposition \ref{surunmodele} shows that $\bar{s}_m$ is consistent as
the maximum likelihood estimator is and that the rate of convergence is
smaller than $\log(T)/T$. It is well known that the MLE is
asymptotically Gaussian in classical settings with a rate of
convergence in $1/T$. But the aim of Proposition \ref{surunmodele} is
not to investigate asymptotic properties: the virtue of the previous
result is its nonasymptotic nature. It allows a dependence of $m$ on
$T$, as soon as (\ref{condgamma}) is satisfied (see Section \ref{sec6}
for the
resulting minimax properties).

There are two terms in the upper bound. The first one $\Vert s_m-s
\Vert^2$
has already been identified as the bias of the projection estimator.
The second term can be viewed as an upper bound for the stochastic or
variance term. Actually, this upper bound is almost sharp. If we assume
that $s$ belongs to $S_m$, that is, $s=s_m$, then the bias disappears
and the quantity $\mathbb{E}(\Vert\bar{s}_m-s \Vert^2)$---a pure
variance term---is in fact upper bounded by a constant times $|m|\log
(T)/T$. But on the
other hand, we have the following result.
\begin{Prop}\label{minimaxun}
Let $m$ be a model such that
$\inf_{I\in m} \ell(I) \geq\ell_0$
then there exists a positive constant $c$ depending on $A,\eta,P,\rho
,H$ such that if $|m|\geq c$ then
\[
\inf_{\hat{s}}\sup_{s\in S_m\cap
\mathcal{L}_{H,P}^{\eta,\rho}}\mathbb{E}_s(\Vert s-\hat{s} \Vert^2)
\geq
\lozenge_{H,P,\eta,\rho,A} \min\biggl(\frac{|m|}{T},\ell_0|m| \biggr).
\]
The infimum over $\hat{s}$ represents the infimum over all the possible
estimators constructed on the observation on $[-A,T]$ of a point
process $(N_t)_t$. $\mathbb{E}_s$ represents the expectation with
respect to
the stationnary Hawkes process $(N_t)_t$ with intensity given by $\Psi_{s}(\cdot)$.
\end{Prop}

Hence, when $s$ belongs to $S_m$, the clipped projection estimator has
a risk which is lower bounded by a constant times $|m|/T$ and upper
bounded by $|m|\log(T)/T$. There is only a loss of a factor $\log(T)$
between the upper
bound and the lower bound. This factor comes from the unboundedness of the
intensity. The best control we can provide for the intensity is to
bound it on
$[0,T]$ by something of the order $\log(T)$. The reader may think to this
really similar fact: the sup of $n$ i.i.d. variables with exponential
moments can only be bounded
with high probability by something of the order $\log(n)$.
Note also that the clipped projection estimator is minimax on $S_m\cap
\mathcal{L}_{H,P}^{\eta,\rho}$ up to this logarithmic term.

Now let us turn to model selection, oracle inequalities and penalty choices.
As before if we know $H,\eta$, and $\rho$, then it is reasonable to
consider the clipped penalized projection estimator, $\bar{s}$ for
theoretical purpose. Recall that the penalized projection estimator
$\tilde{s}=(\tilde{\nu},\tilde{h})$ is given by (\ref{estpen}). Then
the clipped penalized projection estimator, $\bar{s}=(\bar{\nu},\bar
{h})$, is given, for all positive $t$, by
%
%
\begin{equation}\label{clippenalise}
\left\{\begin{array}{l}
\bar{\nu} = \cases{\tilde{\nu}, &\quad if $\rho\leq\tilde{\nu} \leq
\eta$,\cr
\rho, &\quad if $\tilde{\nu} <\rho$,\cr
\eta, &\quad if $\tilde{\nu} >\eta$,}
\\[16pt]
\bar{h} (t) = \cases{\tilde{h}(t), &\quad if $0\leq\tilde{h}(t)\leq H$,
\cr
0, &\quad if $\tilde{h}(t) <0$,\cr
H, &\quad if $\tilde{h}(t) >H$.}
\end{array}\right.
\end{equation}
The next theorem provides
an oracle inequality in expectation [see (\ref{oracleessai})].
\begin{Th}
\label{oracleineqenesp}
Let $(N_t)_{t\in\mathbb{R}}$ be a Hawkes process with intensity $\Psi
_{s}(\cdot)$. Assume
that we know that $s$ belongs to $\mathcal{L}_{H,P}^{\eta,\rho}$.
Moreover,
assume that all the models in $\mathcal{M}_{T}$ are written on $\Gamma
$, a regular
partition of
$(0,A]$ such that (\ref{condgamma}) holds. Let $Q>1$.
Then there exists a positive constant $\kappa$ depending on $\eta
,\rho,P,A,H$
such that
if
%
%
\begin{equation}\label{penth}
\forall m \in\mathcal{M}_{T}\qquad \operatorname{pen}(m)=\kappa
Q(|m|+1)\frac{\log(T)^2}{T},
\end{equation}
then
\begin{eqnarray*}
\mathbb{E}(\Vert\bar{s}-s \Vert)^2 &\leq& \lozenge_{\eta,\rho
,P,A,H}\inf
_{m\in
\mathcal{M}_{T}} \biggl[
\Vert s-s_m \Vert^2 +
(|m|+1)\frac{\log(T)^2}{T} \biggr]\\
&&{} +\lozenge_{\eta,\rho,P,A,H}
\frac{\#\{\mathcal{M}_{T}\}}{T^Q}.
\end{eqnarray*}
\end{Th}

The form of the penalty is a constant times $|m| \log(T)^2/T$, that
is, it
is equal to the variance term up to some logarithmic factor. Remark
also that choosing the penalty as a constant times the dimension leads
to an oracle inequality in expectation. The multiplicative constant is
not an absolute constant but something that depends on all the
parameters that were introduced ($H, \eta, P$, etc.). This is actually
classical. Even in the Gaussian nested case (see \cite{minimal}),
Mallows' $C_p$ multiplicative constant is $2\sigma^2$ where $\sigma^2$
is the variance of the Gaussian noise. The form is simpler than in our
case but still an unknown parameter $\sigma^2$ appears. With respect
to the Gaussian case, remark that there is also some loss due to
logarithmic terms. Finally, for readers who are familiar with model
selection techniques, we do not refine the penalty with the use of
weights, because the concentration formulas we use to derive the
penalty expression are not concentrated enough to allow a real
improvement by using those weights. The Gaussian concentration
inequalities do not apply to Hawkes processes, even if there are some
attempts at proving similar results \cite{RB-Roy}. As a consequence, we
are not able to treat families of models as complex as in
\cite{gauss}.
This lack of concentration actually comes from an obvious essential
feature of the Hawkes' process: its dependency structure. This has
already been noted in several papers on counting processes (see
\cite{reynaudspl} and \cite{reynaudBer}). Here, the dependance is not a
nuisance parameter but the structure we want to estimate via the
function $h$.  Related works may be found in discrete time:
autoregressive process in \cite{bcv} or \cite{bcv2} and Markov chain in
\cite{lacour}. In all these papers, multiplicative constants, which are
usually unknown by practitioners, appear in the penalty term, as in the
Gaussian framework, where the variance noise $\sigma^2$ is usually unknown.
In the Gaussian case, there have been several papers dealing with the
precise theoretical calibration of those constants in a data-driven way
(see \cite{am} or \cite{minimal}). Here, since the concentration
inequalities are too rough, we cannot prove theoretical calibration. So
we have decided to find at least a practical data-driven calibration of
this multiplicative constant (see Section \ref{sec4}).

\section{Practical data-driven calibration via simulations}\label{sec4}
The main drawback of the previous theoretical results is that the
multiplicative constant in the penalty is not computable in practice.
Even if the formula for the factor $\kappa$ is known, it depends heavily
 on the extra knowledge of parameters ($H, \eta, P$, etc.) that
cannot be guessed in practice. On the contrary, $A$ is a meaningful
quantity, at least for our biological purpose. The aim of this section
is to find a performant implementable method of selection, based on the
following theoretical fact: (\ref{penth}) proves that a constant times
the dimension of the model should work.

\subsection{Compared methods}\label{sec41}

Since our simulation design (see Section \ref{sec43}) is computationally
demanding, we restricted ourselves to models $m$ with at most 15
intervals. Consequently, we did not consider the \textit{Nested strategy}
because it would only involve five models in the family. We then only
focus on the three following strategies: \textit{Regular}, \textit{Irregular} and
\textit{Islands}.
Since we are looking for a penalty that is inspired by (\ref{penth}),
we compare our penalized methods to the most naive approach, namely the
Hold-out procedure described below. As stated in the
\hyperref[sec1]{Introduction}, the log-likelihood contrast coupled with an AIC penalty
(see, e.g., \cite{fado}) is only adapted to functions $g$
defined on regular partitions, so we do not consider this method here.
Moreover, the truncated estimators are designed
for minimax theoretical purposes, but of course they depend on parameters
($H$, etc.) that cannot be guessed in practice. They also force the
estimate of $h$ to be nonnegative. Therefore, in this section, we only
use nontruncated estimators [see (\ref{chaps}), (\ref{critpen}), (\ref{estpen})].

\textit{Hold-out}. The naive approach is based on the following fact
(which can
be made completely and theoretically explicit in the self-exciting
case). We know
(see Lemma \ref{contra}) that $\gamma_T$ is a contrast. We know also that
$\mathbb{E}(\gamma_T(f))=\Vert f-s \Vert_D^2-\Vert s \Vert_D^2$.
Moreover, we
know that
the projection
estimators $\hat{s}_m$ behave nicely (see Proposition \ref
{surunmodele}). Now
we would like to select a model $\hat{m}$ such that $\hat{s}_{\hat
{m}}$ is as
good as the best possible $\hat{s}_m$. So one
way to select a good model $m$ should be to observe a second
independent Hawkes process
with the same $s$ and to compute the minimizer of
$\gamma_{T,2}(\hat{s}_m)$ over $\mathcal{M}_{T}$ (where $\hat{s}_m$
is computed
with the
first process and $\gamma_{T,2}$ is our contrast but computed with the
second process). However, we do not have in practice two independent Hawkes
processes at our disposal. But one can cut $[-A,T]$ in two almost independent
pieces. Indeed, the points of the process in $[-A,T/2-A]$ and in $[T/2,T]$
can be equal to those of independent stationary Hawkes processes and
this with high probability (see
\cite{RB-Roy}). Hence, in the sequel whenever the Hold-out estimator is
mentioned, and whatever the family $\mathcal{M}_{T}$ is, it is
referring to the
following procedure.
\begin{enumerate}
\item Cut $[-A,T]$ into two pieces: $H_1$ refers to the points of the
process on $[-A,T/2-A]$, $H_2$ refers to the points of the
process on $[T/2,T]$.
\item Compute $\hat{s}_m$ for all the $m$ in $\mathcal{M}_{T}$ by minimizing
the least-square contrast $\gamma_{T,1}$ on $S_m$ computed with only
the points of
$H_1$, that is,
\[
\forall f \in\mathbb{L}^{2}\qquad \gamma_{T,1}(f)= -\frac{2}{T} \int_{0}^{T/2-A}
\Psi_{f}(t) \,dN_t + \frac{1}{T} \int_{0}^{T/2-A}
\Psi_{f}(t)^2 \,dt.
\]

\item Compute $\gamma_{T,2}(\hat{s}_m)$ where $\gamma_{T,2}$ is
computed with $H_2$, that is,
\[
\forall f \in\mathbb{L}^{2}\qquad \gamma_{T,2}(f)= -\frac{2}{T} \int_{T/2+A}^{T}
\Psi_{f}(t) \,dN_t + \frac{1}{T} \int_{T/2+A}^{T}
\Psi_{f}(t)^2 \,dt
\]
and find $\hat{m}=\arg\min_{m\in\mathcal{M}_{T}} \gamma
_{T,2}(\hat{s}_m)$.
\item The Hold-out estimator is defined by
$\tilde{s}^{\mathrm{HO}}:=\hat{s}_{\hat{m}}$.
\end{enumerate}

\textit{Penalized}. Theorem \ref{oracleineqenesp} shows that theoretically
speaking a penalty of the type $K (|m|+1)$ should
work. However, the theoretical multiplicative constant is not only not
computable, it is also too large for practical
purpose. So one needs to consider Theorem
\ref{oracleineqenesp} as a result that guides our intuition toward
the right shape of penalty and one should not consider it as a sacred
and not
improvable way of penalizing. Therefore, we investigate two ways of
calibrating the multiplicative constants.
\begin{enumerate}
\item The first one follows the conclusions of \cite{minimal}. In the
\textit{Regular
strategy}, there exists at most one model per dimension. If there
exists a true
model $m_0$,
then for $|m|$ large (larger than $|m_0|$) $\gamma_T(\hat{s}_m)$ should
behave like
$-k(|m|+1)$. So there is a ``minimal penalty'' as defined by Birg\'e and
Massart of the form $\operatorname{pen}_{\min}= k(|m|+1)$. In this
situation, their
rule is to take
$\operatorname{pen}(m)=2*\operatorname{pen}_{\min}(m)$.

We find a $\hat{k}$ by doing a least-square regression for large
values of
$|m|$ so that
\[
\gamma_T(\hat{s}_m)\simeq-\hat{k}(|m|+1).
\]
Then we take
\[
\hat{m}=\mathop{\arg\min}_{m\in\mathcal{M}_{T}}\gamma_T(\hat{s}_m)+2\hat
{k}(|m|+1),
\]
and we define $\tilde{s}^{\min}:=\hat{s}_{\hat{m}}$.

Let us remark that the framework of \cite{minimal} is Gaussian and
i.i.d. It is, in our opinion, completely out of reach to extend these
theoretical results here. However, at least in the \textit{Regular
strategy}, the
concentration formula that lies at the heart of our proof is really
close to
the one used in \cite{minimal}, 
which tends to prove
that their method could work here.

For the \textit{Irregular} and \textit{Islands strategy}, as a preliminary
step, we need to find the best data-driven model per
dimension, that is,
\[
\hat{m}_D=\mathop{\arg\min}_{m\in\mathcal{M}_{T}, |m|=D}\gamma_T(\hat{s}_m).
\]
Then one can plot as a function of $D$, $\gamma_T(\hat{s}_{\hat{m}_D})$.
In \cite{minimal}, they also obtain another kind of minimal penalty of
the form
$\operatorname{pen}_{\min}= k(D+1)(\log(|\Gamma|/D)+5)$ when the
\textit{Irregular
strategy} is used. But for very small values of $|\Gamma|$ (as here), we
would not see the difference between this form of penalty and the linear
form. Moreover, theoretically speaking, we are not able to justify,
even heuristically, such a form
of penalty for large values of $|\Gamma|$. Indeed, the concentration
formula in our
case 
is quite different for such a complex family. 

So we have decided that we will use the same penalty as before even in
the \textit{Irregular} and \textit{Islands strategies}. That is to say that
we find a $\hat{k}$ by doing a least-square regression for large value of
$D$ so that
\[
\gamma_T(\hat{s}_{\hat{m}_D})\simeq-\hat{k}(D+1).
\]
Then we take
\[
\hat{m}=\mathop{\arg\min}_{m\in\mathcal{M}_{T}}\gamma_T(\hat{s}_m)+2\hat
{k}(|m|+1),
\]
and we define $\tilde{s}^{\min}:=\hat{s}_{\hat{m}}$ even for the
\textit{Irregular} and \textit{Islands strategies}.

\item On the other hand, the choice of $\hat{m}$ by $\tilde{s}^{\min}$
was not
completely satisfactory when using the \textit{Islands} or \textit
{Irregular strategies} (see
the comments on the simulations hereafter). But on the contrast curve:
$D\to
\gamma_T(\hat{s}_{\hat{m}_D})$, we could see a perfectly clear angle
at the true
dimension. So we have decided to compute
$-\bar{k}=\frac{\gamma_T(\hat{s}_{\Gamma})-\gamma_T(\hat
{s}_{\hat
{m}_1})}{|\Gamma|-1}$
and to choose
\[
\hat{m}=\mathop{\arg\min}_{m\in\mathcal{M}_{T}}\gamma_T(\hat{s}_m)+\bar
{k}(|m|+1).
\]
We define $\tilde{s}^{\mathrm{angle}}:=\hat{s}_{\hat{m}}$.
This seems to be a proper automatic way to obtain this angle without
having to
look at the contrast curve. It is still based on the fact that a
multiple of the dimension
should work. This has only been implemented for the \textit{Irregular}
and \textit{Islands
strategies}.

%
%
\begin{table}[b]
\tablewidth=250pt
\caption{Table of the different methods}\label{methods}
\begin{tabular*}{\tablewidth}{@{\extracolsep{\fill}}lcc@{}}
\hline
\textbf{Methods} & \textbf{Strategy} & \textbf{Selection}\\
\hline
1 & \textit{Regular} $N=15$ & Minimal penalty $\tilde{s}^{\min}$\\
2 & \textit{Irregular} $|\Gamma|=15$ & Angle method
$\tilde{s}^{\mathrm{angle}}$\\
3 & \textit{Irregular} $|\Gamma|=15$ & Minimal penalty
$\tilde{s}^{\min}$\\
4 & \textit{Islands} $|\Gamma|=15$ & Angle method
$\tilde{s}^{\mathrm{angle}}$\\
5 & \textit{Islands} $|\Gamma|=15$ & Minimal penalty
$\tilde{s}^{\min}$\\
6 & \textit{Regular} $N=15$ & Hold-out $\tilde{s}^{\mathrm{HO}}$\\
7 & \textit{Irregular} $|\Gamma|=15$ & Hold-out $\tilde{s}^{\mathrm{HO}}$\\
8 & \textit{Islands} $|\Gamma|=15$ & Hold-out $\tilde{s}^{\mathrm{HO}}$\\
\hline
\end{tabular*}
\end{table}

This angle method may be viewed as the ``extension'' of the $L$-curve
method in
inverse problems where one chooses the tuning parameter at the point
of highest curvature.
\end{enumerate}
Table \ref{methods} summarizes our 8 different estimators.

\subsection{Simulated design}\label{sec42}

We have simulated Hawkes processes
with parameters $(\nu,h)$, with $\nu$ in $\{0.001, 0.002, 0.003, 0.004,
0.005\}$,
$h$ having a bounded support in $(0,1000]$ (i.e., $A=1000$) and on a
sequence of
length $[-A,T]$ with $T=100\mbox{,}000$ or $T=500\mbox{,}000$. The fact that the
process is or
not stationary does not seem to influence our procedure with this relatively
short memory (indeed $T\geq100 A$).

The functions $h$ have been designed so that we can see the influence
of $p$ (\ref{defp}) on the
estimation procedure. So $f_1=0.004\mathbh{1}_{[200,400]}$ is a piecewise
constant nonnegative function on the regular partition $\Gamma$
($|\Gamma|=15$)
with integral $0.8$ and we have tested $h=c*f_1$ with $c$ in $\{0.25,
0.5, 0.75, 1\}$ (i.e., $p=0.2, 0.4, 0.6$ and $0.8$, respectively). We
have also tested a
possibly negative function
$f_2=0.003\mathbh{1}_{[200,800/3]}-0.003\mathbh{1}_{[2000/3,2200/3]}$
that is
piecewise constant on $\Gamma$. Note that (see
Section \ref{sec24}) the sign of $h$ should not affect the method (penalized
least-square criterion) whereas the log-likelihood may have some problems
each time $\Psi_{f}(\cdot)$ remains negative on a large interval. The
parameter of
importance here is the integral of the absolute value, which is here
$\int|f_2|=0.8$ and we have tested $h=f_2$. Finally, the method itself
should not be affected by a smooth function $h$: we
have used $f_3$ a nonnegative continuous function (in fact the mixture
of two Gaussian densities) with
integral equal to $0.8$ and we have tested once again $h=f_3$.

Remark that the mean number of observed points belongs to
$[125,12\mbox{,}500]$ which corresponds to the number of occurrences we could
observe in biological data.

\subsection{Implementation}\label{sec43}
The minimization of $\gamma_T$ is actually quite easy since we use a
least-square contrast.
From a matrix point of view, one can associate to some $f$ in $S_m$
[see (\ref{modele})] a vector of $D+1=|m|+1$ coordinates
\[
{\bolds\theta}_m=\pmatrix{ \mu\cr a_{I_1} \cr\vdots\cr a_{I_D}},
\]
where ${I_1},\ldots,{I_D}$ represent the successive intervals of the
model $m$.
Let us introduce
\[
\mathbf{b}_m=\pmatrix{ \displaystyle\frac{1}{T} N_{[0,T]} \vspace*{2pt}\cr
\displaystyle\frac{1}{T} \int_0^T
\Psi
_{(0,\mathbh{1}_{I_1})}(t) \,dN_t \vspace*{2pt}\cr\vdots
\vspace*{2pt}\cr
\displaystyle\frac{1}{T} \int_0^T \Psi_{(0,\mathbh{1}_{I_D})}(t) \,dN_t}
\]
and
\begin{eqnarray*}
\mathbf{X}_m&=&
{\fontsize{8}{10}\selectfont{\left(\matrix{
1 & \frac{1}{T} \int_0^T \Psi_{(0,\mathbh{1}_{I_1})}(t) \,dt & \cdots & \frac{1}{T} \int_0^T \Psi_{(0,\mathbh{1}_{I_D})}(t) \,dt
\cr
\frac{1}{T} \int_0^T \Psi_{(0,\mathbh{1}_{I_1})}(t) \,dt & \frac{1}{T} \int_0^T \Psi_{(0,\mathbh{1}_{I_1})}^2(t) \,dt & \cdots
& \frac{1}{T} \int_0^T \Psi_{(0,\mathbh{1}_{I_1})}(t)\Psi_{(0,\mathbh{1}_{I_D})}(t)
\,dt\cr
\vdots & \vdots &\ddots  & \vdots\cr
\frac{1}{T} \int_0^T \Psi_{(0,\mathbh{1}_{I_D})}(t) \,dt & \frac{1}{T}
\int_0^T \Psi_{(0,\mathbh{1}_{I_1})}(t)\Psi_{(0,\mathbh{1}_{I_D})}(t) \,dt &
\cdots & \frac{1}{T} \int_0^T \Psi_{(0,\mathbh{1}_{I_D})}^2(t)
\,dt}\right)}}.
\end{eqnarray*}
It is not difficult to see that the contrast $\gamma_T(f)$ can be
written
\[
\gamma_T(f)=-2{\bolds\theta}_m \mathbf{b}_m + {^{t}{\bolds\theta}_m}
\mathbf{X}_m
{\bolds\theta}_m.
\]
Therefore, the minimizer $\hat{\bolds\theta}_m$ of $\gamma_T(f)$ over
$f$ in $S_m$ satisfies $\mathbf{X}_m \hat{\bolds\theta}_m=\mathbf{b}_m$,
that is, $\hat{\bolds\theta}_m=\mathbf{X}_m^{-1} \mathbf{b}_m$.
Since the functions $\Psi_{(0,\mathbh{1}_{I})}(t)$ are piecewise
constants, despite their randomness, it may be long but not that
difficult to compute $\mathbf{X}_m$. It is also possible to compute
$\mathbf{X}_\Gamma$ and to deduce from it the different $\mathbf
{X}_m$'s, when one
uses the Islands or Irregular strategies.
Nevertheless, both \textit{Islands} and \textit{Irregular strategies} require
to calculate each vector $\hat{\bolds\theta}_m$ for the $2^{|\Gamma
|}$ possible models $m$ and to store them to evaluate the oracle risk
(see below).
We thus restricted our Monte Carlo simulations to models $m$ with less
than 15 intervals. For the analysis of single real data sets, the
technical limitation of our programs is $|\Gamma|=26$ due to the
$2^{|\Gamma|}$ possible models. The programs have been implemented in
\texttt{R} and are available upon request.

\subsection{Results}\label{sec44}

The quality of the estimation procedures is measured thanks to two
criteria: the risk of the estimators and the associated oracle ratio.
\begin{itemize}
\item
We call \textit{Risk} of an estimator the Mean Square Error of this estimator
over 100
simulations, that is, we compute for each simulation $\Vert s-\hat{s}
\Vert^2$ and next we
compute the average over 100 simulations. Note that with the range of our
parameters, the error of estimation of $\nu$ will be really negligible
with respect to the error of estimation for $h$, so that $\Vert s-\hat
{s} \Vert^2\simeq\int_0^A (h-\hat{h})^2$.
\item
The \textit{Oracle Risk} is for each method the minimal risk, that is,\break $\min
_{m\in\mathcal{M}_{T}} \mathit{Risk}(\hat{s}_m)$. All our methods give an estimator
$\tilde
{s}$ that is selected among a
family of $ \hat{s}_m$'s. The \textit{Oracle Ratio} is the ratio of the risk of
$\tilde{s}$ divided by the Oracle
Risk, that is,
\[
\frac{\mathit{Risk}(\tilde{s})}{\min_{m\in\mathcal{M}_{T}} \mathit{Risk}(\hat{s}_m)}.
\]
If the \textit{Oracle Ratio} is $1$, then the risk of $\tilde{s}$ is the one of
the best estimator in the family.
Note that the definition of $\mathcal{M}_{T}$ and even the definition
of $\hat
{s}_m$ appearing in the \textit{Oracle Ratio}
may change from one method to another one.
\end{itemize}

\begin{figure}

\includegraphics{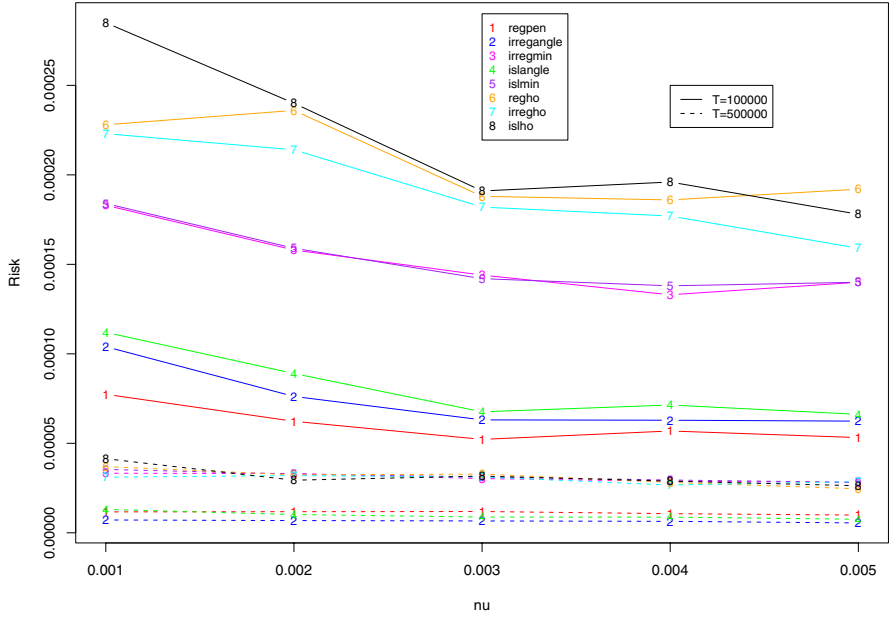}

\caption{Risk of the 8 different methods for $h=0.5*f_1$ for different
values of $\nu$ and $T$.}\label{riskf1c05}
\end{figure}

\begin{figure}

\includegraphics{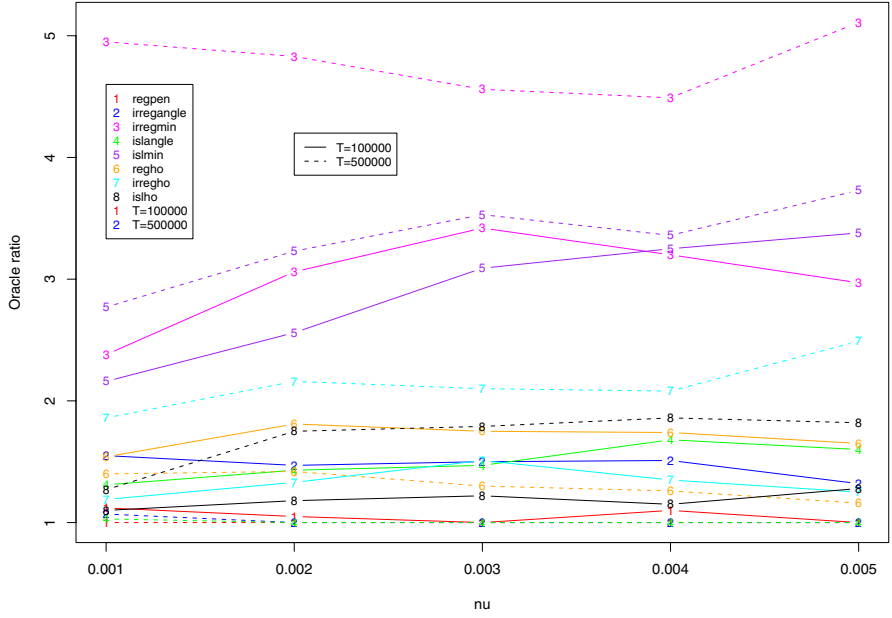}

\caption{Oracle Ratio of the 8 different methods for $h=0.5*f_1$ for
different values of $\nu$ and $T$.}\label{oraclef1c05}
\end{figure}

Figure \ref{riskf1c05} gives the \textit{Risk} of our estimators for
$h=0.5*f_1$ for
various $\nu$ and $T$. We first clearly see
%
%
that the risk decreases when $T$ increases whatever the method. Then we
see that the
``best methods'' are methods 1, 2 and 4, that is, the \textit{Regular
strategy} with minimal penalty
and the \textit{Irregular} and \textit{Islands strategies} with the
angle method.
For the \textit{Irregular} and \textit{Islands strategies}, the minimal
penalty seems to behave like the Hold-out strategies.
There seems also to be a slight improvement when $\nu$ becomes larger, tending
to prove that, if the mean total number of points $\mathbb
{E}(N[0,T])=\nu
T/(1-p)$ grows, the estimation is improved---at least in our range of
parameters. Figure \ref{oraclef1c05} gives the \textit{Oracle Ratio} of our estimators
in the same context. The \textit{Oracle Ratio} is really close to 1 for methods
1, 2 and 4 when
$T=500\mbox{,}000$ whatever $\nu$ is.
Remark that the \textit{Oracle Ratio} for the Hold-out
estimators (methods 6, 7 and 8) is not that large, but since the
estimators $\hat{s}_m$ are computed
with half of the data, their \textit{Risks} are not as small as the projection
estimators used
in the penalty methods. This explains why the \textit{Risk} of the Hold-out
methods is large when the \textit{Oracle Ratio} is close to 1. The \textit{Oracle
Ratio}
is improving when $T$
becomes larger for our three favorite methods (namely 1, 2, 4).

%
%
%

%
\begin{figure}

\includegraphics{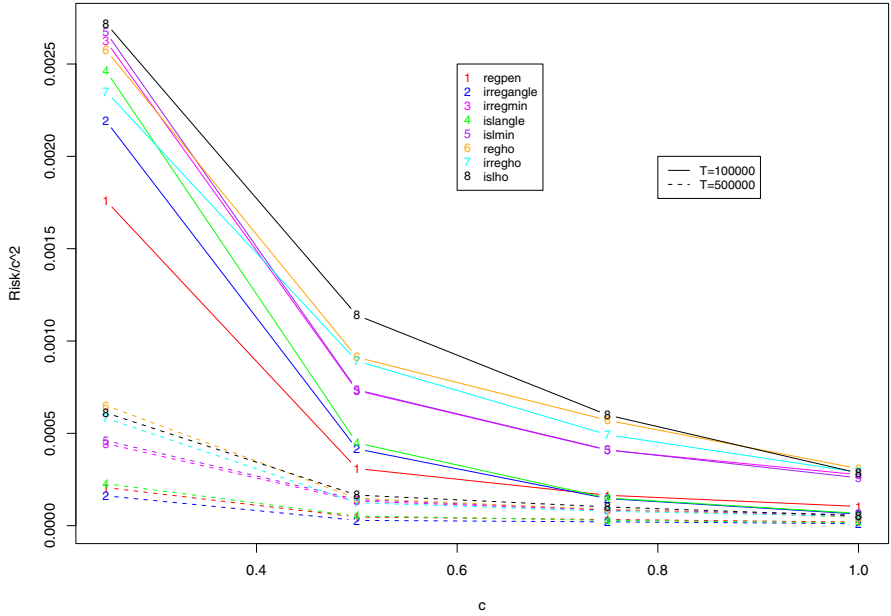}

\caption{Rescaled Risk ($\mathit{Risk}/c^2$) of the 8 different methods for
$h=c*f_1$ and $\nu=0.001$, for different values of $c$ and $T$.}
\label{riskscalef1nu0001}
\end{figure}

%
\begin{figure}

\includegraphics{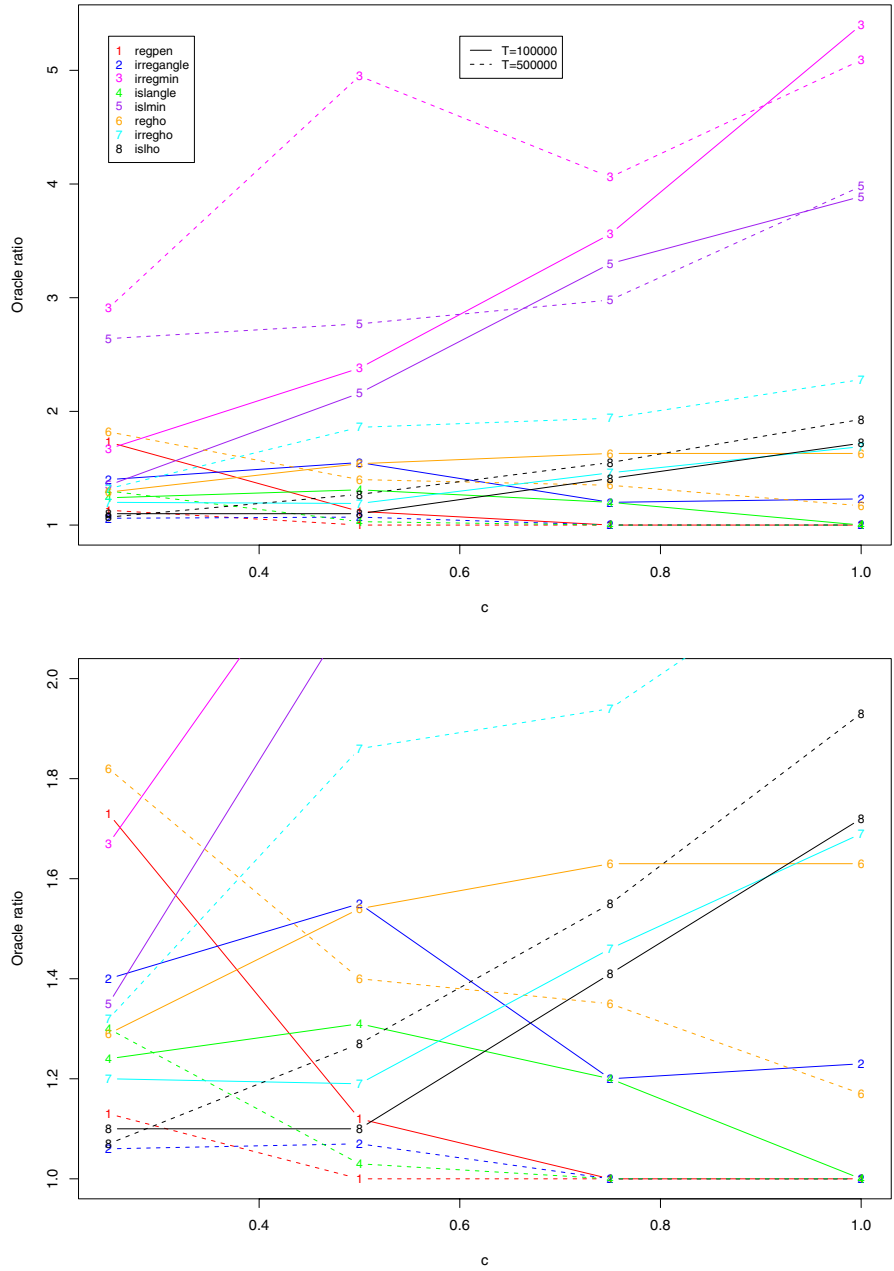}

\caption{Oracle Ratio of our estimators for $h=c*f_1$ and $\nu=0.001$
for different values of $c$ and $T$ (top). The bottom picture zooms in
on the top picture for Oracle Ratio between 1 and 2.}\label{oraclef1nu0001}
\end{figure}

%
\begin{figure}

\includegraphics{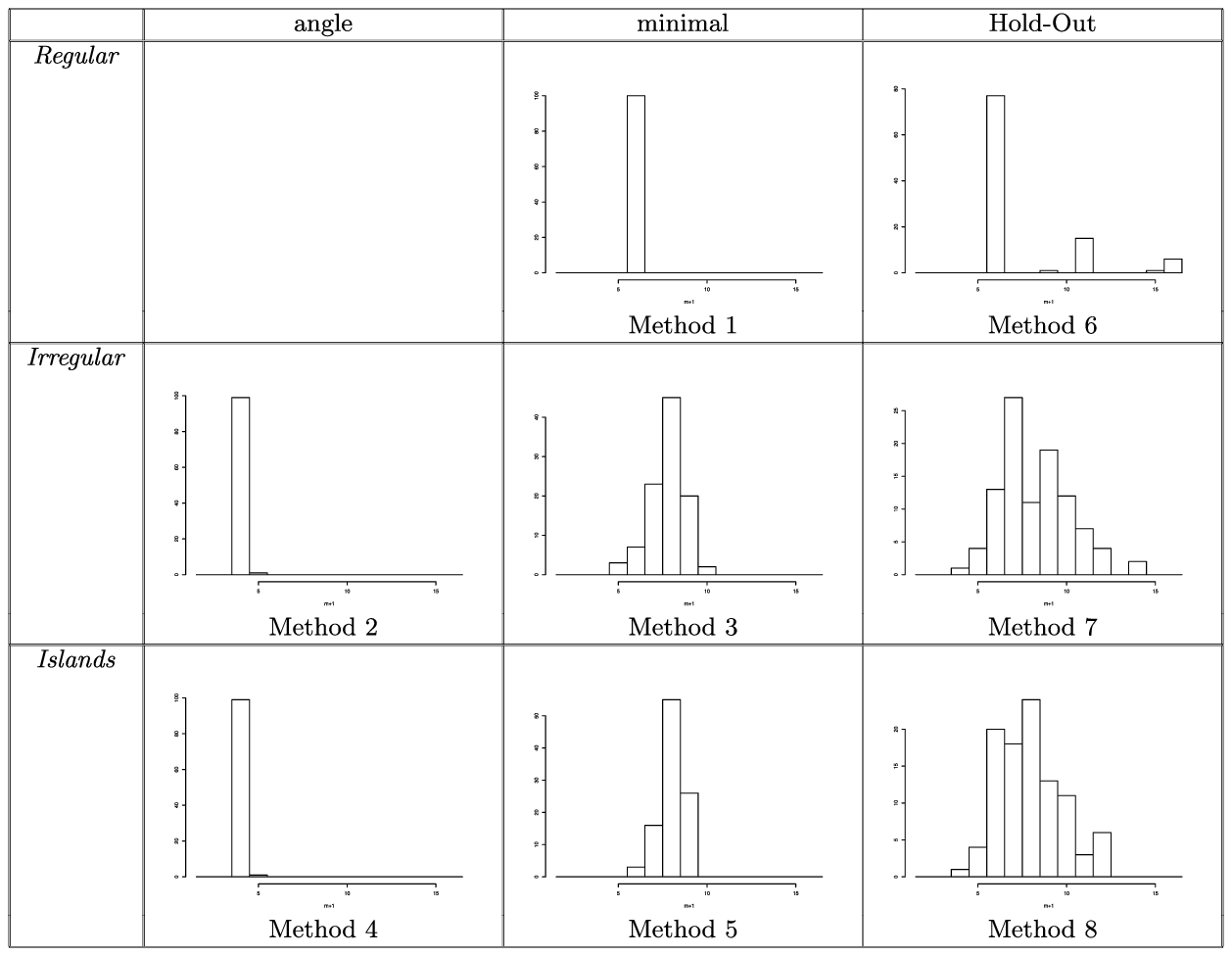}

\caption{Frequency of the chosen dimension $|\hat{m}|+1$ for the different
methods when $T=500\mbox{,}000$, $\nu=0.001$ and $h=0.5*f_1$. Note that the
true dimension is 6 for the Regular method (chosen
in 100\% of the simulations by method 1) and 4 for the
Irregular and Islands methods (chosen in more than
95\% of the simulations by methods 2 and 4).}
\label{histogrammedim}
\end{figure}

Figure \ref{riskscalef1nu0001} gives the variation of the \textit{Risk} with
respect to $p$ (\ref{defp}). Since $h=c*f_1$ and since $c$ varies, the
\textit{Rescaled} \textit{Risk}, $\mathit{Risk}/c^2$, gives (up to some negligible term
corresponding to $\nu$) the risk of $\tilde{h}/c$ as an estimator of
$f_1$. We
clearly see that when $T$ or $c$ becomes larger the \textit{Rescaled} \textit{Risk} is
decreasing. So it definitely seems that if the mean total number of
points grows, the estimation is improving. Methods 1, 2 and 4 seem to
be still the
more precise ones. Figure \ref{oraclef1nu0001} gives the
\textit{Oracle Ratio} in the same situation. Once again there is an improvement when
$T$ grows at least for our three favorite methods (1, 2 and 4) and the \textit{Oracle
Ratio} is 1 when $T=500\mbox{,}000$ and $c=0.8$. The same comment
about a good \textit{Oracle Ratio} for the Hold-out methods apply.

Figure \ref{histogrammedim} gives the frequency of the chosen dimension,
namely $|\hat{m}|+1$ for the different methods. Clearly, methods 1, 2
and 4 are
correctly choosing the true dimension in most of the simulations when the
other methods overestimate the true dimension.

%
\begin{figure}

\includegraphics{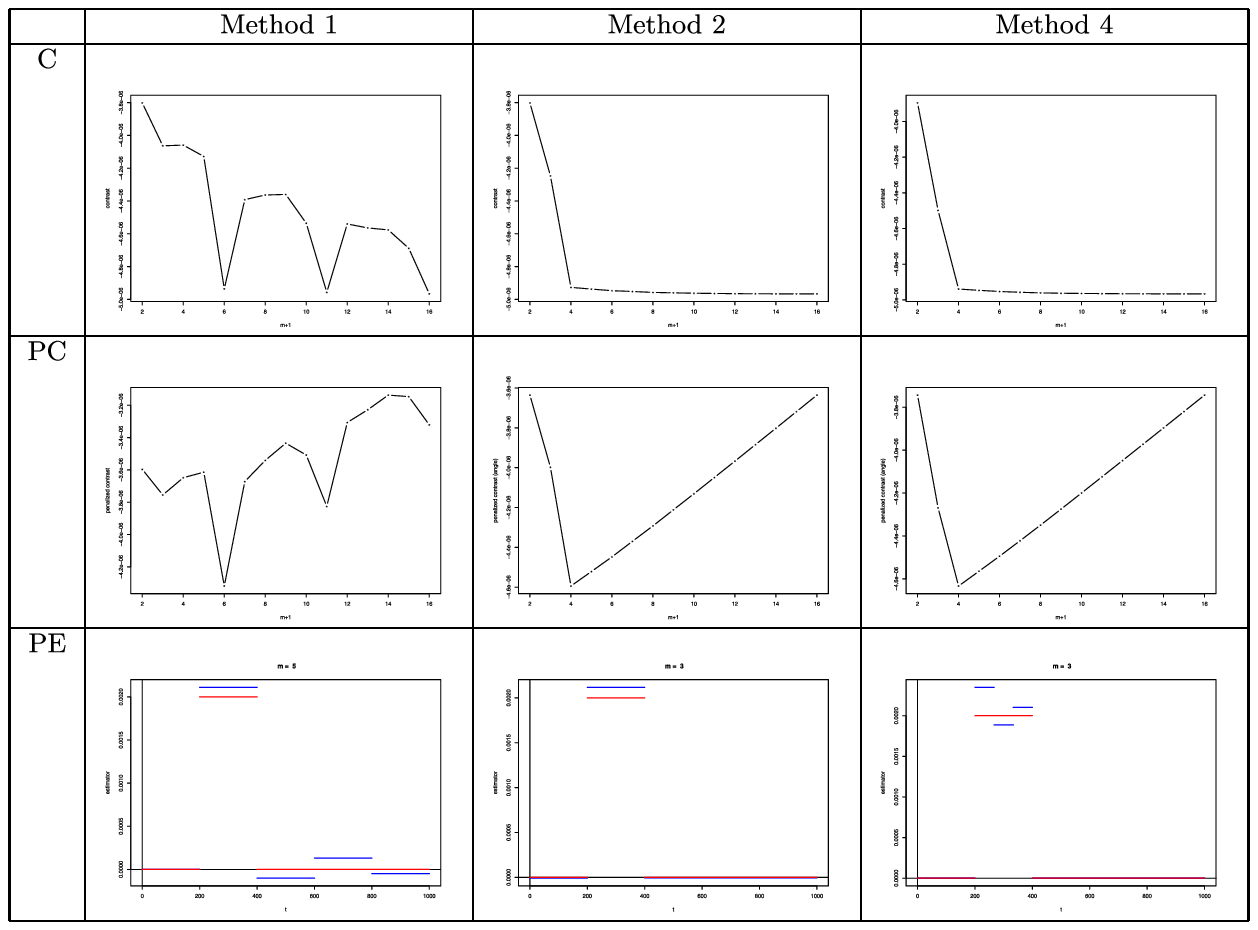}

\caption{Contrast (C) and penalized contrast (PC) as a function of the
dimension for
the three favorite methods on one simulation with $T=500\mbox{,}000$, $\nu=0.001$
and $h=0.5*f_1$. The chosen estimators (PE) are in blue whereas the
function $h=0.5*f_1$ is in red.}
\label{contrastefig}
\end{figure}

Finally, Figure \ref{contrastefig} shows the resulting estimators of
methods 1,
2 and 4 on one simulation. In particular, before penalizing, note that one
clearly sees an angle on the contrast curve at the true dimension and that
penalizing by the angle method (methods 2 and 4) gives an automatic way to
find the position of this angle.

%
\begin{figure}

\includegraphics{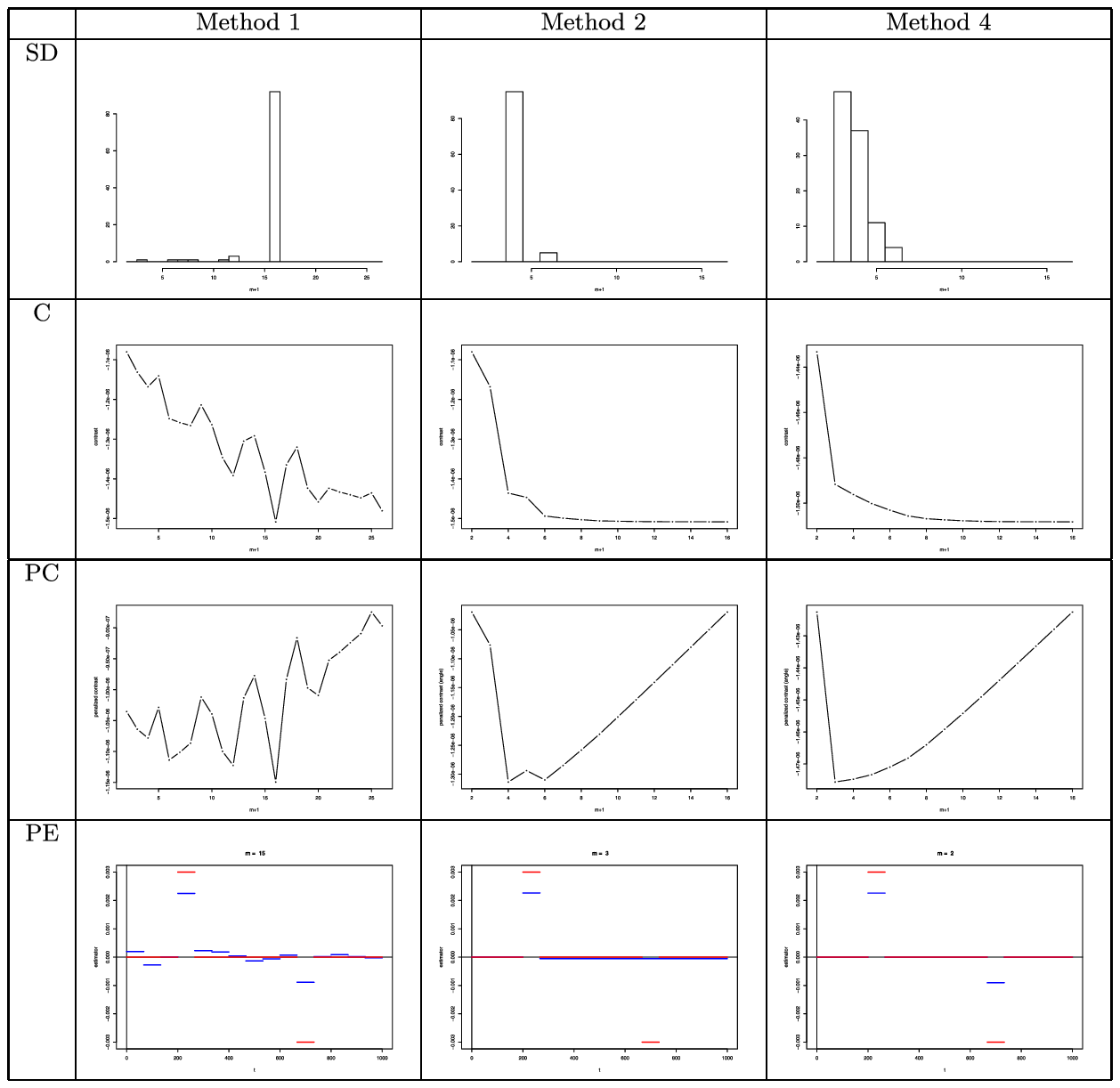}

\caption{Histogram of the selected dimension over 100 simulations
(SD). Contrast (C) and penalized contrast (PC) as a function of the
dimension for
the three favorite methods on one simulation with $T=500\mbox{,}000$, $\nu=0.001$
and $h=f_2$. The true dimension is 16 for method 1
(Regular), 6 for method 2
(Irregular) and 3 for method 4 (Islands). The chosen
estimators (PE) are in blue whereas the
function $h=f_2$ is in red.}\label{contrastefigf2}
\end{figure}

Figure \ref{contrastefigf2} shows the results for the possibly
negative function $f_2$ and only for our three favorite methods
(1, 2, 4). For this function only, and because the true dimension is 16
for method 1, we use for method 1, $|\Gamma|=25$. Note that (i)
methods~1 and~4 select the right dimension whereas method 2
(\textit{Irregular strategy}) does not see the negative jump and that
(ii) it is also more easy to detect the precise position of
the fluctuations on the sparse estimate given by method 4
(compared to method~1). For sake of simplicity, we do not give
the \textit{Risk} values, but it is sufficient to note that, for all the
methods, they are small (with a slight advantage for method~4) and that
the \textit{Oracle Ratios} are close to 1.

%
\begin{figure}

\includegraphics{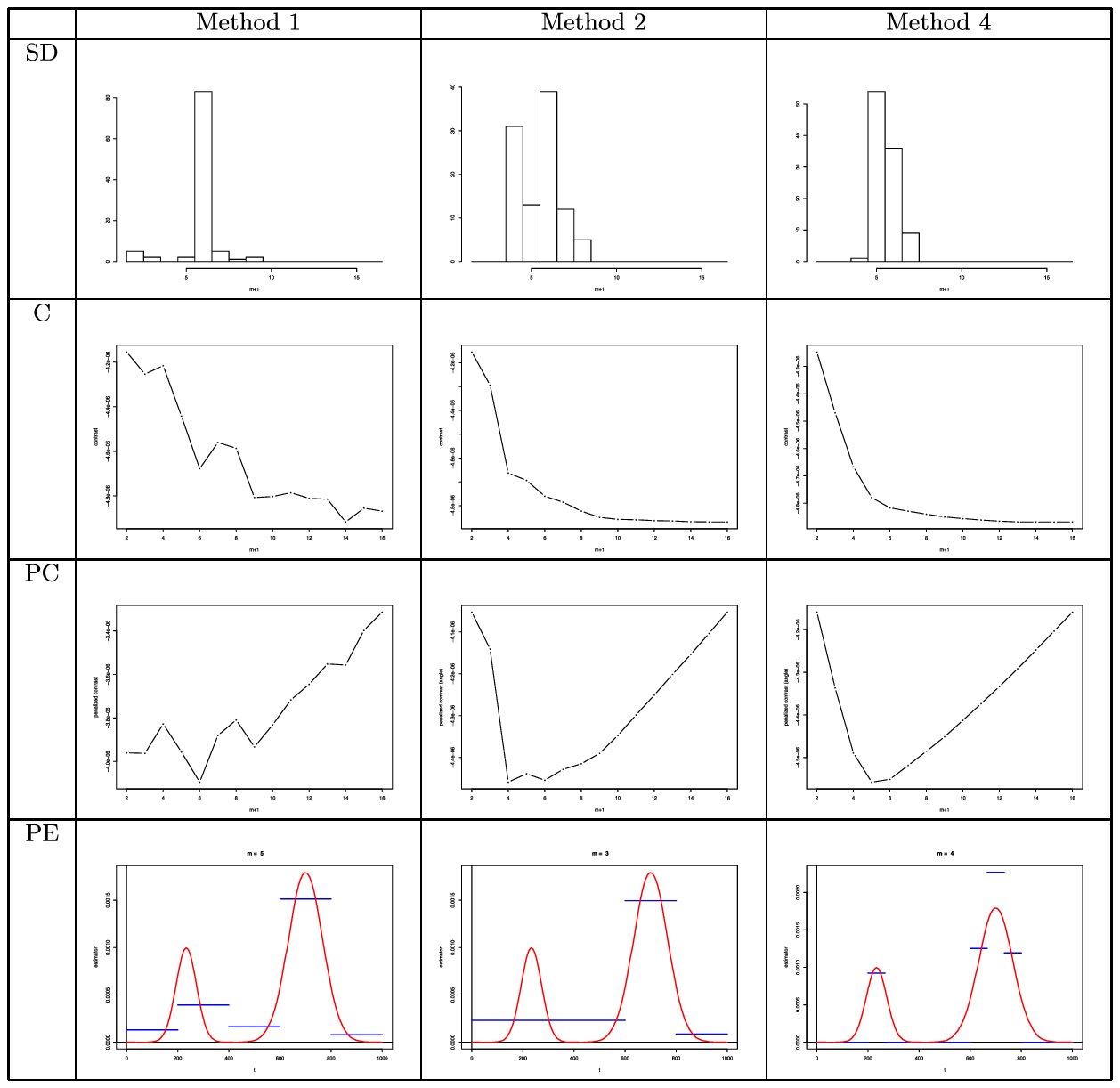}

\caption{Histogram of the selected dimension over 100 simulations (SD).
Contrast (C) and penalized contrast (PC) as a function of the dimension for
the three favorite methods on one simulation with $T=500\mbox{,}000$, $\nu=0.001$
and $h=f_3$. The chosen estimators (PE) are in blue whereas the
function $h=f_3$ is in red.}\label{contrastefigf3}
\end{figure}

Figure \ref{contrastefigf3} gives the same results for the smooth
function $f_3$. Of course, since the projection estimators are piecewise
constant, they cannot look really close to $f_3$. But in any case,
method 1 and more interestingly method 4 gives the right position for
the spikes whereas method 2 does not see the smallest bump.

Finally, let us conclude the simulations by noting that the penalized
projection estimators with the \textit{Islands strategy} and the angle
penalty (method 4) seems to be an appropriate method for detecting
local spikes and bumps in the function $h$ and even negative jumps.

\section{Applications on real data}\label{sec5}

We have applied the penalized (angle\break method) estimation procedure with
the \textit{Island strategy} (method 4) to two data sets related to
occurrences of genes
or DNA motifs along both strands of the complete genome of the
bacterium \textit{Escherichia coli} ($T=9\mbox{,}288\mbox{,}442$). In both cases, we
used $A=10\mbox{,}000$ as the longest dependence between events and
the finest partition corresponds to $|\Gamma|=15$.

\begin{figure}

\includegraphics{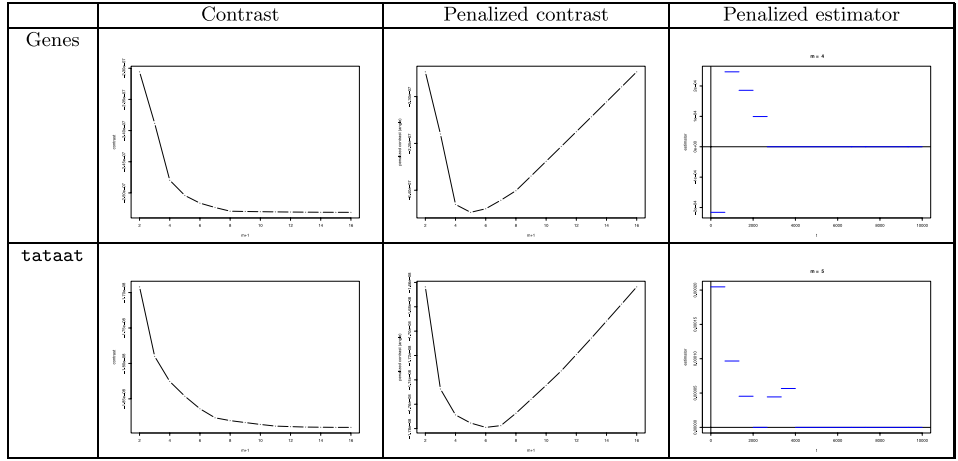}

\caption{Contrasts, penalized contrasts and chosen estimators for both
E. coli datasets.}
\label{applifig}
\end{figure}

\begin{figure}[b]

\includegraphics{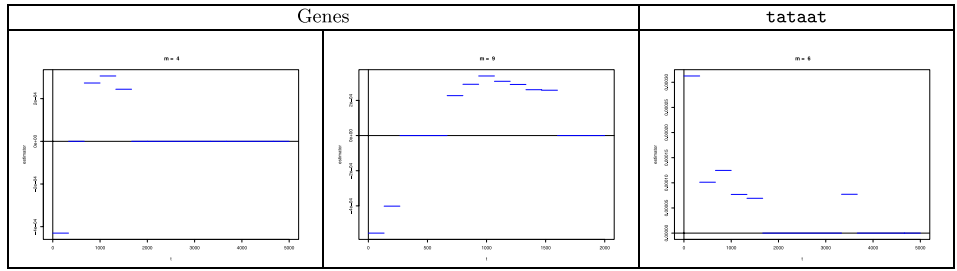}

\caption{Chosen estimators for both E. coli datasets for
different values of $A$: $A=5000$ (left, right) and $A=2000$ (middle).}
\label{applifig2}
\end{figure}

The first process corresponds to the occurrences of the 4290 genes.
Figure \ref{applifig} (top) gives the associated contrast and penalized
contrast, together with the chosen estimator of $h$ ($\hat{m}=4$
and $\hat{\nu}=3.64$ 10$^{-4}$). The shape of this estimator tells
us that:
\begin{itemize}
\item gene occurrences seem to be uncorrelated down to 2600 basepairs,
\item they are avoided at a short distance ($\sim$0--$500$ bps) and
\item favored at distances $\sim$700--$2000$ bps apart.
\end{itemize}
This general trend has been refined by shortening the support
$A$ to 5000 and then to 2000 (see Figure \ref{applifig2}). It then
clearly appears both a negative effect at distances less than 250 bps,
and a positive one around 1000 bps. This
is completely coherent with biological observations: genes on the
same strand do not usually overlap, they are about 1000 bps long in average,
and there are few intergenic regions along bacterial genomes (compact
genomes).

%

The second process corresponds to the 1036 occurrences of the DNA
motif \texttt{tataat}. Figure \ref{applifig} (bottom) gives the
associated contrast and penalized
contrast, together with the chosen estimator of $h$ ($\hat{m}=5$
and $\hat{\nu}=7.82$ $10^{-5}$). The shape of the estimator suggests
that:
\begin{itemize}
\item occurrences seem to be uncorrelated down to 4000 basepairs,
\item favored at distances $\sim0$--$1500$ bps and 3000 bps apart,
\item highly favored at a short distance apart (less than 600 bps).
\end{itemize}
After shortening the support $A$ to 5000 (see Figure
\ref{applifig2}), the shape of the chosen estimator shows that there
actually are 3 types of favored distances: very short distances (less
than 300
bps), around 1000 bps and around 3500 bps.
This trend is again coherent with the fact that
(i) the motif \texttt{tataat} is self-overlapping (two successive
occurrences can occur at a distance 5 apart),
(ii) this motif is part of the most common
promoter of \textit{E. coli} meaning that it should occur in front of
the majority of the genes (and these genes seem to be favored at
distances around 1000~bps apart from the previous example),
(iii) some particular successive genes (operons) can
be regulated by the same promoter (this could explain the third bump).

%

%
\begin{figure}

\includegraphics{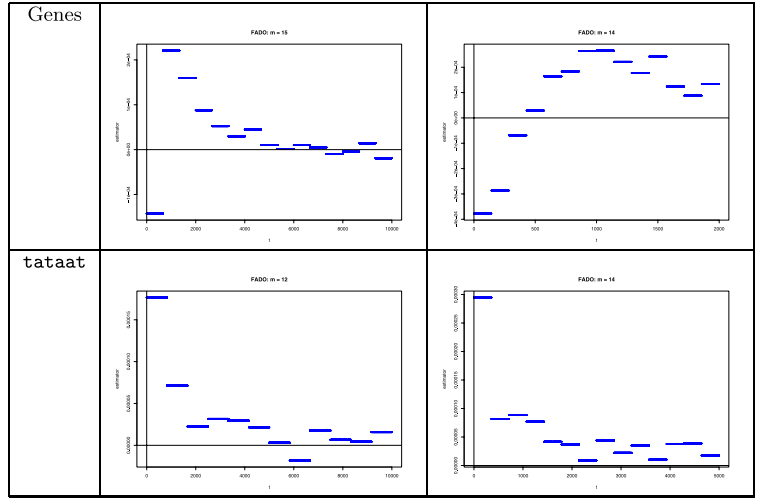}

\caption{FADO estimators for both E. coli datasets for
different values of $A$: $A=10\mbox{,}000$ (left) and $A=2000$ (right) for
genes or $A=5000$ (right) for \texttt{tataat}.}
\label{avecfado}
\end{figure}

Figure \ref{avecfado} presents the results of the FADO procedure
\cite{fado}. Here, we have forced the estimators to be piecewise constant to
make the comparison easier. Note, however, that the FADO procedure may be
implemented with splines of any fixed degree.

Our results are in agreement with the ones obtained by
FADO.
 Our new approach has two advantages. First, it
gives a better idea of the support $A$ of the function $h$: indeed, the
estimator provided by FADO (cf. Figure \ref{avecfado} top-left) has some
fluctuations until the end of the interval whereas our estimator (cf.
Figure \ref{applifig} top-right) points out that nothing significant
happens after 3000 bps. Second, our method leads
to models of smaller dimension ($|m|=4$ for \textit{Islands} versus
$|m|=15$ for FADO). The limitation of our method is
essentially that we only consider piecewise constant estimators, but
this is enough to get a general trend on favored or avoided distances
within a point process.

\section{Minimax properties}\label{sec6}
The theoretical procedures of Proposition \ref{surunmodele} and Theorem
\ref{oracleineqenesp} have more theoretical properties than just an
oracle inequality. This section provides their minimax properties. In
particular, even if it has not been implemented for technical reasons
that were described above, the \textit{Nested strategy} leads to an
adaptive minimax estimator. Such kinds of estimators were not known in
the Hawkes model, as far as we know.

\subsection{\texorpdfstring{H\"olderian functions}{Holderian functions}}\label{sec61}
First, one can prove the following lower bound.
\begin{Prop}\label{lowerholder}
Let $L>0$ and $1\geq a >0$. Let
\[
\mathcal{H}_{L,a}=\{s=(\nu,h)\in\mathbb{L}^{2} / \forall x,y \in(0,A],
|h(x)-h(y)|\leq
L|x-y|^a\}.
\]
Then
\[
\inf_{\hat{s}}\sup_{s\in\mathcal{H}_{L,a}\cap
\mathcal{L}_{H,P}^{\eta,\rho}}\mathbb{E}_s(\Vert s-\hat{s} \Vert^2)
\geq
\lozenge_{H,P,A,\eta,\rho,a}\min\bigl(L^{{2}/({2a+1})} T^{-
{2a}/({2a+1})},1 \bigr).
\]
The infimum over $\hat{s}$ represents the infimum over all the possible
estimators constructed on the observation on $[-A,T]$ of a point
process $(N_t)_t$. $\mathbb{E}_s$ represents the expectation with
respect to
the stationnary Hawkes process $(N_t)$ with intensity given by $\Psi_{s}(\cdot)$.
\end{Prop}

But on the other hand, let us consider the clipped projection estimator
$\bar{s}_m$ with $m$ a regular partition of $(0,A]$ such that
\[
|m|\simeq\bigl(T/\log(T)\bigr)^{1/(2a+1)}.
\]
If the function $h$ is in $\mathcal{H}_{L,a}\cap
\mathcal{L}_{H,P}^{\eta,\rho}$ with $a\in(1/2,1]$, then, applying
Proposition~\ref{surunmodele}, $\bar{s}_m$ satisfies
\[
\mathbb{E}(\Vert\bar{s}_m-s \Vert^2)\leq\lozenge_{H,P,A,\eta,\rho
,L,a} \biggl(\frac
{\log
(T)}{T} \biggr)^{2a/(2a+1)}.
\]
Compared with the lower bound of the minimax risk (Proposition
\ref{lowerholder}), we only lose a logarithmic factor:
the clipped projection estimators are minimax on $\mathcal
{H}_{L,a}\cap
\mathcal{L}_{H,P}^{\eta,\rho}$, with $a\in(1/2,1]$, up to some
logarithmic term.
We cannot go beyond $a=1/2$
because one needs $|m|\ll\sqrt{T}$ in Proposition \ref{surunmodele}.

Of course, we need to know $a$ to find $\bar{s}_m$, so $\bar{s}_m$ is
not adaptive with respect to $a$. But the clipped penalized projection
estimator $\bar{s}$ with the \textit{Nested strategy} can be adaptive
with respect to $a$. It is sufficient to take $J\simeq\log_2(\sqrt
{T}/\log(T)^3)$ to guarantee (\ref{condgamma}). Then we apply Theorem
\ref{oracleineqenesp} with $Q=1.1$, for instance. Since $\#\{\mathcal
{M}_{T}\}$ is
of the order $\log(T)$, we obtain that
\[
\mathbb{E}(\Vert\bar{s}-s \Vert)^2 \leq\lozenge_{H,\eta,P,A,\rho}
\inf_{m\in
\mathcal{M}_{T}} \biggl[
\Vert s-s_m \Vert^2 +
(|m|+1)\frac{\log(T)^2}{T} \biggr].
\]
If $h$ is in $\mathcal{H}_{L,a}\cap
\mathcal{L}_{H,P}^{\eta,\rho}$ with $a\in(1/2,1]$, then there exists
$m$ in $\mathcal{M}_{T}$ such that
\[
|m|\simeq\bigl(T/\log(T)^2\bigr)^{1/(2a+1)}
\]
and consequently
\[
\mathbb{E}(\Vert\bar{s}-s \Vert)^2 \leq\lozenge_{H,\eta,P,\rho,A,L,a}
\biggl(\frac
{\log
(T)^2}{T} \biggr)^{2a/(2a+1)}.
\]
Therefore, the clipped penalized projection estimator $\bar{s}$ with
the \textit{Nested strategy} and the theoretical penalty given by (\ref
{penth}) is adaptive minimax on $\{\mathcal{H}_{L,a}\cap
\mathcal{L}_{H,P}^{\eta,\rho}, a\in(1/2,1]\}$ up to some
logarithmic term.

\subsection{Irregular and Islands sets}\label{sec62}
Let us apply Theorem \ref{oracleineqenesp} to the \textit{Irregular
strategy} and \textit{Islands strategy}.
In both cases,
the limiting factor here is $\#\{\mathcal{M}_{T}\}$. Take $N\leq\log
_2(T)$, then
$\#\{\mathcal{M}_{T}\}\leq T$ and if
$Q\geq2$ we obtain that
\[
\mathbb{E}(\Vert\bar{s}-s \Vert)^2 \leq\lozenge_{H,P,A,\eta,\rho}
\inf_{m\in
\mathcal{M}_{T}} \biggl[
\Vert s-s_m \Vert^2 +
(|m|+1)\frac{\log(T)^2}{T} \biggr].
\]
To measure performances of those estimators, one needs to introduce a
set of sparse functions $h$, functions that are difficult to estimate
with a \textit{Nested strategy}. A piecewise function $h$ is usually
thought as sparse if the resulting partition is irregular with few intervals.
So we define the Irregular set by
%
%
\begin{equation}\label{irreg}
S^{\mathrm{irr}}_{\Gamma,D}:=\bigcup_{m \ \mathrm{partition}\ \mathrm{written}\
\mathrm{on}\ \Gamma,\ |m| =D}S_m.
\end{equation}
Then, if $s$ belongs to $S^{\mathrm{irr}}_{\Gamma,D}$, using the \textit{Irregular
strategy}, the clipped penalized projection estimator satisfies
\[
\mathbb{E}(\Vert\bar{s}-s \Vert)^2 \leq\lozenge_{H,P,\eta,\rho,A}
D\frac{\log(T)^2}{T}.
\]
But for our biological purpose, the sparsity lies in the support of
$h$. So we define the Islands set by
%
%
\begin{equation}\label{isla}
S^{\mathrm{isl}}_{\Gamma,D}: =\bigcup_{m\subset\Gamma, |m| =D}S_m.
\end{equation}
Then, if $s$ belongs to $S^{\mathrm{isl}}_{\Gamma,D}$, using the \textit{Islands
strategy}, the clipped penalized projection estimator also satisfies
\[
\mathbb{E}(\Vert\bar{s}-s \Vert)^2 \leq\lozenge_{H,P,\eta,\rho,A}
D\frac{\log(T)^2}{T}.
\]
On the other hand it is possible to compute lower bounds for the
minimax risk over those sets.
\begin{Prop}\label{minimaxislirr}
Let $\Gamma$ be a partition of $(0,A]$ such that
\mbox{$\inf_{I\in\Gamma} \ell(I) \geq\ell_0$}.
Let $|\Gamma|=N$ and let $D$ be a
positive integer such that $N\geq4D$.
If $D\geq c_2(A,\eta,P,\rho,H)>1$, for $c_2$ some positive constant
depending on $A,\eta,P,\rho,H$, then
\[
\inf_{\hat{s}}\sup_{s\in S^{\mathrm{isl}}_{\Gamma,D}\cap
\mathcal{L}_{H,P}^{\eta,\rho}}\mathbb{E}_s(\Vert s-\hat{s} \Vert^2)
\geq
\lozenge_{H,P,A,\eta,\rho}\min\biggl(\frac{D\log({N}/{D}
)}{T},D\ell_0 \biggr)
\]
and
\[
\inf_{\hat{s}}\sup_{s\in S^{\mathrm{irr}}_{\Gamma,D}\cap
\mathcal{L}_{H,P}^{\eta,\rho}}\mathbb{E}_s(\Vert s-\hat{s} \Vert^2)
\geq
\lozenge_{H,P,A,\eta,\rho}\min\biggl(\frac{D\log({N}/{D}
)}{T},D\ell_0 \biggr).
\]
The infimum over $\hat{s}$ represents the infimum over all the possible
estimators constructed on the observation on $[-A,T]$ of a point
process $(N_t)_t$. $\mathbb{E}_s$ represents the expectation with
respect to
the stationary Hawkes process $(N_t)$ with intensity given by $\Psi_{s}(\cdot)$.
\end{Prop}

To clarify the situation, it is better to take $N=|\Gamma|\simeq\log
(T)$. If $D\simeq\log(T)^a$ with $a<1$ then the lower bound on the
minimax risk is of the order $\log(T)^a \log\log T /T$ when the risk of
the clipped penalized projection estimator (for both strategies) is
upper bounded by $\log(T)^{a+2}/T$, and this whatever $a$ is. So our
estimator matches the rate $1/T$ up to a logarithmic term. Of course
the most fundamental part is this logarithmic term. Think, however, that
there exists some function $h$ in those sets, such that the function
belongs to $S_\Gamma$ but to none of the other spaces $S_m$ for $m$ in
the family $\mathcal{M}_{T}$ described by the \textit{Nested strategy}.
Consequently,
a clipped penalized estimator with the \textit{Nested strategy} would
have an upper bound on the risk of the order $\log(T)^3/T$ by applying
Theorem \ref{oracleineqenesp}. So the \textit{Irregular} and \textit
{Islands} strategies have not only good practical properties, but there
is also definitely a theoretical improvement in the upper bound of the risk.

\section{Technical results}\label{sec7}
\subsection{Oracle inequality in probability}\label{sec71}
The following result is actually the one at the origin of Theorem \ref
{oracleineqenesp}.
Note that this result holds for the practical estimator, $\tilde{s}$,
which is not clipped.
\begin{Th}
\label{oracleineqenproba}
Let $(N_t)_{t\in\mathbb{R}}$ be a Hawkes process with intensity $\Psi
_{s}(\cdot)$.
Let $H$,
$\eta$ and $A$ be positive known constants such that $s=(\nu,h)$ satisfies
$\nu\in[0,\eta]$ and $h(\cdot)\in[0,H]$.

Moreover, assume that the family $\mathcal{M}_{T}$ satisfies
\[
\inf_{m\in\mathcal{M}_{T}}\inf_{I\in m}\ell(I)\geq\ell_0>0.
\]

Let
$\mathcal{S}$ be a finite vectorial subspace of $\mathbb{L}^{2}$ containing
all the piecewise constant functions constructed on the models of
$\mathcal{M}_{T}$.
Let $R>r>0$ be positive real numbers, let $\mathcal{N}$
be a positive integer and let us consider the following event:
\[
\mathcal{B}= \bigl\{\forall t \in[0,T], N\bigl([t-A,t)\bigr)\leq
\mathcal{N} \mbox{ and } \forall f\in\mathcal{S}, r^2\Vert f \Vert
^2\leq
D_T^2(f)\leq R^2 \Vert f \Vert^2 \bigr\},
\]
where $N([t-A,t))$ represents the number of points of the Hawkes
process $(N_t)_t$ in the interval $[t-A,t)$.
We set $\Lambda=(\eta+H\mathcal{N})R^2/r^2$ and we consider
$\varepsilon$ and
$x$ any arbitrary positive constants. If
for all $m\in\mathcal{M}_{T}$
\[
\operatorname{pen}(m)\geq(1+\varepsilon)^3\Lambda
\frac
{|m|+1}{T}\bigl(1+3\sqrt{2x}\bigr)^2,
\]
then there exists an event $\Omega_x$ with probability larger than
$1-3\#\{\mathcal{M}_{T}\} e^{-x}$ such that
for
all $m\in\mathcal{M}_{T}$, both following inequalities hold:
%
%
\begin{eqnarray}
\label{BESOIN}\hspace*{28pt}
&&\frac{\varepsilon r^2}{1+\varepsilon}\Vert\tilde{s}-s \Vert^2 \mathbh
{1}_{\mathcal{B}\cap\Omega
_x} \nonumber\\
&&\qquad
\leq(1+\varepsilon) D^2_T(s_m-s) + (1+\varepsilon^{-1})
D_T^2(s-s_\perp)+ r^2\Vert s_\perp-s \Vert
^2\\
&&\qquad\quad{} +\frac{r^2}{1+\varepsilon}
\Vert s-s_m \Vert^2 + (1+\varepsilon)\operatorname{pen}(m)
+ \lozenge_\varepsilon\frac{\Lambda}{T} x+\lozenge
_\varepsilon\frac{1+\mathcal{N}^2/\ell
_0}{r^2T^2}x^2,\nonumber
\end{eqnarray}
where $s_\perp$ denotes the orthogonal projection for $\Vert\cdot \Vert$
of $s$
on $\mathcal{S}$, and
%
%
\begin{eqnarray}
\label{BESOIN2}
&&
r^2\frac{\varepsilon}{1+\varepsilon}\mathbb{E}
(\Vert\tilde{s}-s \Vert^2\mathbh{1}_{\mathcal{B}\cap\Omega
_x} ) \nonumber\\
&&\qquad
\leq\biggl((2+\varepsilon+\varepsilon^{-1})
K^2+\frac{2+\varepsilon}{1+\varepsilon}r^2 \biggr) \Vert s-s_m \Vert^2\\
&&\qquad\quad{} + (1+\varepsilon
)\operatorname{pen}(m)
+\lozenge_\varepsilon\Lambda\frac{x}{T}+ \lozenge
_\varepsilon\frac{1+\mathcal{N}^2/\ell
_0}{r^2T^2}x^2,\nonumber
\end{eqnarray}
where $K$ is a positive constant depending on $s$ such that $\Vert f
\Vert_D\leq K\Vert f \Vert$ for all $f$ in $\mathbb{L}^{2}$ (see
Lemma \ref{norme}).
\end{Th}
\begin{remark}\label{remark1}
This result is really the most fundamental to understand
how the Hawkes
process can be easily handled once we only focus on a nice event, namely
$\mathcal{B}$. We have ``hidden'' in $\mathcal{B}$ the fact that the intensity
of the
process is unbounded: on $\mathcal{B}$, the number of points per
interval of
length $A$
is controlled, so the intensity is bounded on this event. We have also
``hidden'' in $\mathcal{B}$ the fact that we are working with a
natural norm, namely
$D_T$, which is random and which may eventually behave badly: on
$\mathcal{B}$, $D_T$
is equivalent to the deterministic norm $\Vert\cdot \Vert$ for functions in
$\mathcal{S}$.
More precisely, the result of (\ref{BESOIN}) mixes $\Vert\cdot \Vert$ and
$D_T(\cdot)$ but holds in probability. On the contrary, (\ref{BESOIN2}) is
weaker but more readable since it holds in expectation with only one
norm $\Vert\cdot \Vert$. Note also that $\mathcal{B}$ is
observable, so if one observes that we are on $\mathcal{B}$, (\ref{BESOIN2})
shows that a
penalty of the type a factor times the dimension can work really well to
select the right dimension. Indeed, note that if, in the family
$\mathcal{M}_{T}$,
there is a
``true'' model $m$ (meaning that $s=s_m$) and if the penalty is
correctly chosen, then (\ref{BESOIN2})
proves that $\Vert\tilde{s}-s \Vert^2$ is of the same order as the lower
bound on
the minimax risk on $m$, namely $|m|/T$ (see Proposition \ref
{minimaxun} for
the precise lower bound). In that sense, this is an oracle inequality. The
procedure is adaptive because it can select the right model without
knowing it. But of
course this hides something of importance. If $\mathcal{B}$ is not that
frequent, then
the result is completely useless from a theoretical point of view since one
cannot guarantee that the risk of the penalized estimator and
even the risk of the projection estimators themselves are small.
\end{remark}
\begin{remark}\label{remark2}
In fact, we will see in the next subsection that the
choices of
$\mathcal{N},R,r,\mathcal{M}_{T}$ are really important to control
$\mathcal{B}$. In
particular, we are not able
at the end to manage families of models with a very high complexity as in
\cite{gauss} or in most of the other works in model selection (see
Theorem \ref{oracleineqenesp} and Section \ref{sec6}). This is probably
due to a lack of independency and boundedness in the process itself.
\end{remark}
\begin{remark}\label{remark3}
Note also that the oracle inequality in probability
(\ref{BESOIN}) of Theorem \ref{oracleineqenproba}
remains true for the more general process defined by (\ref{sineg}) once
we replace $\mathcal{B}$ by $\mathcal{B}\cap\mathcal{B}'$ where
$\mathcal{B}'=\{\forall t\leq T,
\lambda(t)>0\}$. But of course then, $\mathcal{B}'$ is not
observable. This tends to
prove that even in case of self-inhibition a penalty of the type a
constant times the dimension is working.
\end{remark}

\subsection{Control of $\mathcal{B}$}\label{sec72}

The assumptions of Theorem \ref{oracleineqenesp} are in fact a direct
consequence of the assumptions needed to control $\mathcal{B}$, as
shown in the following result.
\begin{Prop}
\label{controldeb}
Let $s\in
\mathcal{L}_{H,P}^{\eta,\rho}$ and $R$ and $r$ such that
\[
R^2> 2\max\biggl(1, \frac{\eta}{(1-P)^2}\bigl(\eta A+(1-P)^{-1}\bigr) \biggr)
\quad\mbox{and}\quad
r^2< \min\biggl(\frac{\rho}{4}, \frac{1-P}{8A \eta+1} \biggr).
\]
Moreover let
\[
\mathcal{N}=\frac{6 \log(T)}{P-\log P-1}.
\]
Let us finally assume that $\mathcal{S}$, defined in
Theorem \ref{oracleineqenproba}, is included in $S_{\Gamma}$ where
$\Gamma$
is a regular partition of $(0,A]$ such that
\[
|\Gamma|\leq\frac{\sqrt{T}}{(\log T)^3}.
\]
Then, under the assumptions of Theorem \ref{oracleineqenproba}, there
exists $T_0>0$ depending on $\eta,\rho,P,A,R$ and $r$, such that
for all $T>T_0$,
\[
\mathbb{P}(\mathcal{B}^c)\leq\lozenge_{\eta,P,A}\frac{1}{T^2}.
\]
\end{Prop}

These technical results imply very easily Proposition \ref{surunmodele}
and Theorem \ref{oracleineqenesp}.
\begin{pf*}{Proof of Theorem \ref{oracleineqenesp}}
We apply (\ref{BESOIN2}) of Theorem
\ref{oracleineqenproba} to $\tilde{s}$. Since $\bar{s}$ is closer to
$s$ than $\tilde{s}$,
the inequality is also true for $\bar{s}$. We choose $x=Q\log(T)$ and
$\mathcal{N}$,
$R, r$ according to Proposition \ref{controldeb}.
On the complement of $\mathcal{B}\cap\Omega_x$, we bound $\Vert
\bar{s}-s \Vert$ by $\eta^2+H^2A$
and the probability of the complement of the event by
\[
\lozenge_{\eta,P,A,\rho,H}
\biggl(\frac{1}{T^2}
+\frac{{\#\{\mathcal{M}_{T}\}}}{T^Q}
\biggr).
\]
The same control may be applied if $T$ is not large enough. To complete
the proof, note finally that $K\leq\lozenge_{\eta,P,A}$.
\end{pf*}
\begin{pf*}{Proof of Proposition \ref{surunmodele}}
We can apply Theorem \ref{oracleineqenproba} to a family that is
reduced to only one
model $m$. If the inequality is true for the nontruncated estimator, and
if we
know the bounds on $s$ then the inequality is necessarily true for the
truncated estimator, which is closer to $s$ than $\tilde{s}$. Then the
penalty is not needed to compute the estimator but it
appears nevertheless in both oracle inequalities. We can conclude by
similar arguments as Theorem \ref{oracleineqenesp}, but if we take
$x=\log(T)$ in~(\ref{BESOIN2}), we lose a logarithmic factor with
respect to Proposition \ref{surunmodele}. We actually obtain
Proposition \ref{surunmodele} by integrating also in $x$ the oracle
inequality in probability (\ref{BESOIN}) and we conclude by similar
arguments, using that $\Vert\cdot \Vert_D\leq K
\Vert\cdot \Vert$.
\end{pf*}

\section{Sketch of proofs for 
the technical and minimax results}\label{sec8}
\subsection{Contrast and norm}\label{sec81}
First, let us begin with a result that makes clear the link between the
classical properties of the Hawkes process (namely the Bartlett
spectrum) and the quantity $\int g^2$ that is appearing in the
definition of the $\mathbb{L}^{2}$ space (\ref{l2}).
\begin{Lemma}
\label{calculE2}
Let $(N_t)_{t\in\mathbb{R}}$ be a Hawkes process with intensity $\Psi
_{s}(\cdot)$.
Let $g$ be a function on $\mathbb{R}_+$ such
that $\int_0^{+\infty} g(u) \,du$ is finite. Then for all $t$,
\begin{eqnarray*}
&&\mathbb{E}\biggl[ \biggl(\int_{-\infty}^tg(t-u) \,dN_u \biggr)^2 \biggr]
\\
&&\qquad= \frac{\nu^2}{(1-p)^2}
\biggl(\int
_0^{+\infty} g(u) \,du \biggr)^2 + \int_{\mathbb{R}}
|\mathcal{F}{{g}}(-w)|^2 f_N(w) \,dw
\\
&&\qquad\leq\frac{\nu^2}{(1-p)^2} \biggl(\int_0^{+\infty} g(u) \,du \biggr)^2 + \frac
{\nu
}{(1-p)^3} \int_0^{+\infty} g^2(u) \,du,
\end{eqnarray*}
where
\[
f_N(w) = \frac{\nu}{2 \pi(1 - p) |1 -
\mathcal{F}{{h}}(w)|^2}
\]
is the spectral density of $(N_t)_{t\in\mathbb{R}}$.
\end{Lemma}
\begin{remark*}[(Notation)]
$\mathcal{F}{{h}}$ is the Fourier transform of $h$, that is, $\mathcal
{F}{{h}}(x) =
\int_\mathbb{R} e^{ixt} h(t) \,dt$.
\end{remark*}
\begin{pf*}{Proof of Lemma \ref{calculE2}}
Let $\phi_t(u) = \mathbh{1}_{u <t} g(t-u)$. We know (see
\cite{BrM01}, page 123) that
\[
\operatorname{Var}\biggl[\int_{\mathbb{R}} \phi_t(u) \,dN_u \biggr] = \int
_{\mathbb{R}} |\mathcal{F}{{\phi}}_t(w)|^2
f_N(w)\, dw.
\]
Moreover, since $g$ has a positive support,
$
\mathcal{F}{{\phi}}_t(w) = e^{iwt} \mathcal{F}{{g}}(-w).
$
Hence,
\[
\operatorname{Var}\biggl[ \int_{\mathbb{R}} \phi_t(u) \,dN_u \biggr] = \int
_{\mathbb{R}}
|\mathcal{F}{{g}}(-w)|^2 f_N(w) \,dw.
\]
But we also know that (see \cite{HaO74})
\[
\lambda=
\mathbb{E}(\lambda(t))=\frac{\nu}{1-p}.
\]
Consequently,
\begin{eqnarray*}
\mathbb{E}\biggl[
\biggl(\int_{-\infty}^{t}g(t-u)\,dN_u \biggr)^2 \biggr] &=& \operatorname{Var}\biggl[
\int_{\mathbb{R}}\phi_t(u) \,dN_u \biggr] + \biggl( \mathbb{E}\biggl(\int_{\mathbb
{R}} \phi_t(u) \,dN_u \biggr) \biggr)^2
\\
&=& \operatorname{Var}\biggl[ \int_{\mathbb{R}}\phi_t(u) \,dN_u \biggr]+
\biggl(\lambda\int_{0}^{+\infty} g(u) \,du \biggr)^2,
\end{eqnarray*}
which gives the first part of the lemma. The second part is due to
Plancherel's identity, which states
%
%
\begin{equation}
\label{Plancherel}
\int_{\mathbb{R}}|\mathcal{F}{{g}}(-w)|^2 \, dw= 2\pi\int_0^A g^2(x)\,dx,
\end{equation}
and the fact that $f_N$ is upper bounded by
$\nu/[2\pi(1-p)^3]$ since $h$ is nonnegative.
\end{pf*}

Lemma \ref{calculE2} is at the root of Lemma \ref{norme}, which gives
the equivalence between the $\mathbb{L}^{2}$-norms, $\Vert\cdot \Vert$ and
$\Vert\cdot \Vert_D$,
equivalence that is essential for our analysis. Lemma \ref{calculE2}
essentially represents the main feature of the lengthy but necessary
computations of Lemma \ref{norme}.
The proof of Lemma \ref{norme} is consequently omitted and can be found
in~\cite{sopacomplet}.
\begin{Lemma}
\label{norme}
The functional $D^2_T$ is a quadratic form on
$\mathbb{L}^{2}$ and its expectation $\Vert\cdot \Vert_D^2$ [see (\ref
{normD})] is the
square of a norm on $\mathbb{L}^{2}$
satisfying
%
%
\begin{equation}\label{eqvnorm}
\forall f \in\mathbb{L}^{2}\qquad L \Vert f \Vert\leq\Vert f \Vert_D
\leq K \Vert f \Vert,
\end{equation}
where
\[
K^2= 2\max\biggl[1, \frac{\nu}{(1-p)^2} \biggl(\nu A+\frac{1}{1-p} \biggr) \biggr] \quad\mbox{and}
\quad L^2 = \min\biggl[\frac\nu4, \frac{1-p}{8A \nu+1} \biggr].
\]
\end{Lemma}

Lemma \ref{norme} has a direct corollary: $\gamma_T$ defines a contrast.
\begin{Lemma}
\label{contra}
Let $(N_t)_{t\in\mathbb{R}}$ be a Hawkes process with intensity $\Psi
_{s}(\cdot)$.
Then the
functional given by
\[
\forall f \in\mathbb{L}^{2}\qquad \gamma_T(f)= -\frac{2}{T} \int
_{0}^{T} \Psi_{f}(t)
\,dN_t + \frac{1}{T} \int_{0}^{T}
\Psi_{f}(t)^2 \,dt
\]
is a contrast, that is, $\mathbb{E}(\gamma_T(f))$ is minimal for $f=s$.
\end{Lemma}
\begin{pf}
Let us compute $\mathbb{E}(\gamma_T(f))$. As $\lambda(t)=\Psi
_{s}(t)$, one can write
by the martingale properties of $dN_t-\Psi_{s}(t)\,dt$ using the associate
bilinear form of $D^2_T(f)$
that
\begin{eqnarray*}
\mathbb{E}(\gamma_T(f))&=& \mathbb{E}\biggl[-\frac{2}{T}\int_0^T \Psi_{f}(t)
\,dN_t \biggr]+\mathbb{E}(D^2_T(f))\\
&=& \mathbb{E}\biggl[-\frac{2}{T}\int_0^T \Psi_{f}(t) \Psi_{s}(t) \,dt
\biggr]+\Vert f \Vert_D^2\\
&=& \Vert f-s \Vert_D^2-\Vert s \Vert_D^2.
\end{eqnarray*}
Consequently, $\mathbb{E}(\gamma_T(f))$ is minimal when $f=s$ since Lemma
\ref{norme} proves that $\Vert\cdot \Vert_D$ is a norm.
\end{pf}

\subsection[Proof of Theorem 2]{Proof of Theorem \protect\ref{oracleineqenproba}}\label{sec82}
This proof is quite classical in model selection. It heavily depends on
a concentration inequality for $\chi^2$-type statistics that has been
derived in \cite{reynaudspl} and which holds for any counting process.
The main feature is to use the martingale properties of $ N_t-\int_0^t
\lambda(u) \,du$ [see (\ref{class})]. We do not need any further
properties of the Hawkes process to obtain (\ref{BESOIN}) (see Remark
\ref{remark3}).

We give here a sketch of the proof to emphasize that:
\begin{enumerate}
\item the oracle inequalities of Theorem \ref{oracleineqenproba} hold
for $\tilde{s}$ the
practical estimator and not only the clipped one, and
\item that (\ref{BESOIN}) holds for possible negative function $h$ up
to a minor correction (see Remark \ref{remark4} at the end of the proof).
\end{enumerate}
More details may be found in \cite{sopacomplet}.
\begin{pf*}{Proof of Theorem \ref{oracleineqenproba}}
Let $m$ be a fixed partition of $\mathcal{M}_{T}$.
By construction, we obtain
%
%
\begin{equation}
\label{majorationcontraste}
\gamma_T(\tilde{s}) + \operatorname{pen}(\hat{m})\leq\gamma
_T(\hat{s}_m) +
\operatorname{pen}(m)
\leq\gamma_T(s_m) + \operatorname{pen}(m).
\end{equation}
Let us denote for all $f$ in $\mathbb{L}^{2}$,
\[
\nu_T(f)= \frac{1}{T} \int_0^T \Psi_{f}(t) \bigl(dN_t-\Psi_{s}(t)\,dt\bigr),
\]
which is linear in $f$.
Then (\ref{contraste}) becomes $\gamma_T(f)= D^2_T(f-s)-D^2_T(s)-2
\nu_T(f)$ and (\ref{majorationcontraste}) leads to
%
%
\begin{equation}
\label{depart}
D^2_T(\tilde{s}-s) \leq D^2_T(s_m-s) + 2 \nu_T(\tilde{s}-s_m) +
\operatorname{pen}(m)
- \operatorname{pen}(\hat{m}).
\end{equation}
By linearity of $\nu_T$, $\nu_T(\tilde{s}-s_m) = \nu_T(\tilde
{s}-s_{\hat{m}})+
\nu_T(s_{\hat{m}}-s_m)$. Now let us control each term in the right-hand
side of (\ref{depart}).
\begin{enumerate}
\item Let us begin with $A_1=2\nu_T(\tilde{s}-s_{\hat{m}})$. For all
$m'$ in $\mathcal{M}_{T}$, we set
%
%
\begin{equation}
\label{w}
W_{m'} = \sup_{f\in S_{m'}} \frac{\nu_T(f)}{\Vert f \Vert}.
\end{equation}
Thus, $A_1 \leq2\Vert\tilde{s}-s_{\hat{m}} \Vert W_{\hat{m}}$.
Therefore, for all $\theta>0$, one has the following upper bound:
%
%
\begin{equation}
\label{mona}
A_1 \leq\theta\Vert\tilde{s}-s_{\hat{m}} \Vert^2+ \frac
{1}{\theta}
W_{\hat{m}}^2.
\end{equation}
Now we need to control $W_{\hat{m}}$ which is doubly random: for fixed
$m$, $W_m$ is
random but the choice $\hat{m}$ is random too. So one needs to
control each $W_{m'}$'s to control $W_{\hat{m}}$.

To do so, we first need to find
a simpler form for $W_{m'}$. Note that
\[
\{(1,0) \}\cup\biggl\{ \biggl(0, \frac{\mathbh{1}_I}{\sqrt{\ell(I)}} \biggr), I\in
m' \biggr\}
\]
is an orthonormal basis of $S_{m'}$ for $\Vert\cdot \Vert$. For all $I\in
m'$, let us denote
\[
N_I(t)=\Psi_{(0,\mathbh{1}_I)}(t).
\]
Then
we can prove that (see \cite{sopacomplet})
\[
\hspace*{-4pt}W_{m'}=\sqrt{ \Biggl( \int_0^T \frac{1}{T}
\bigl(dN_t-\Psi_{s}(t)\,dt\bigr) \Biggr)^2 +\sum_{I\in m'} \Biggl( \int_0^T\frac
{N_I(t)}{T\sqrt
{\ell(I)}}\bigl(dN_t-\Psi_{s}(t)\,dt\bigr) \Biggr)^2}.
\]
Let $\mathcal{T}$ be defined by
\[
\mathcal{T}= \biggl\{t\geq0 / N\bigl([t-A,t)\bigr)>\mathcal{N} \mbox{ or } \exists f
\in\mathcal{S},
\frac1T\int_0^t\Psi_{f}(u)^2 \,du > R^2 \Vert f \Vert^2 \biggr\}
\]
and let
$\tau$ be the stopping time defined by
\[
\tau=\inf\{t\geq0, t\in\mathcal{T}\}.
\]
It is quite easy to see that if $t$ belongs to $\mathcal{T}$ then
there exists
$t'<t$ such that $t'$ belongs to $\mathcal{T}$. Hence, $\tau$ does
not belong
to $\mathcal{T}$ and since $\int_0^t\Psi_{f}(u)^2 \,du$ is increasing
in $t$,
saying that we restrict ourselves to $\mathcal{B}$ implies that $\tau
\geq T$. Finally,
we can write that on $\mathcal{B}$,
$W_{m'}=Z_{m'}$ defined by
\begin{eqnarray*}
Z_{m'} &=& \biggl( \biggl( \int_0^T \frac{1}{T} \mathbh{1}_{t\leq\tau}
\bigl(dN_t-\Psi_{s}(t)\,dt\bigr) \biggr)^2\\
&&\hspace*{5.8pt}{} +\sum_{I\in m'} \biggl(
\int_0^T\frac{N_I(t)}{T\sqrt{\ell(I)}}\mathbh{1}_{t\leq
\tau}\bigl(dN_t-\Psi_{s}(t)\,dt\bigr) \biggr)^2\biggr)^{1/2}.
\end{eqnarray*}
Written in this way, this is a $\chi^2$-type statistics as defined in
\cite{reynaudspl}, since the $N_I(\cdot)$'s are predictable processes and
so is
$\mathbh{1}_{t\leq\tau}$. So Corollary 2 of \cite{reynaudspl} gives
that
with probability larger than $1-2e^{-x}$,
\[
Z_{m'} \leq\sqrt{C_{m'}} + 3 \sqrt{2 v x} +b x,
\]
where
\[
C_{m'}= \int_{0}^{T} \biggl[\frac{1}{T^2} + \sum_{I\in m'}
\frac{N_I^2(t)}{ T^2 \ell(I)} \biggr] \mathbh{1}_{t\leq\tau}\Psi_{s}(t)
\,dt,\qquad v= \Vert C_{m'} \Vert_{\infty},
\]
and where $b$ is a deterministic constant that should satisfy
\[
b^2 \geq\mathbh{1}_{t\leq\tau} \biggl[\frac{1}{T^2} +
\sum_{I \in m'} \frac{N_I^2(t)}{ T^2 \ell(I)} \biggr].
\]
%
Once we are restricted to $\{\tau\geq T\}$, we can use the quantities
defined in $\mathcal{B}$ to upper bound $C_{m'}$, $v$ and $b$ (see
details in \cite{sopacomplet}). Finally, on $\mathcal{B}$, with probability
larger than $1- 2\#\{\mathcal{M}_{T}\}e^{-x}$,
%
%
\begin{equation}
\label{majow}
W_{\hat{m}}\leq\sqrt{(\eta+ H \mathcal{N})
R^2\frac{|\hat{m}|+1}{T}} \bigl(1+3\sqrt{2x} \bigr)
+ \frac{\sqrt{1+\mathcal{N}^2/\ell_0}}{T}x.
\end{equation}
Let us fix some positive numbers $\theta$ and $\varepsilon$ that will
be chosen later and let us go back to $A_1$. We obtain the following
upper bound:
%
%
\begin{eqnarray}
\label{majoa}
A_1 &\leq& \theta\Vert\tilde{s}-s_{\hat{m}} \Vert^2+ \frac
{1}{\theta}
\biggl[ (1+\varepsilon) (\eta+ H \mathcal{N})
R^2\frac{|\hat{m}|+1}{T} \bigl(1+3\sqrt{2x}
\bigr)^2\nonumber\\[-8pt]\\[-8pt]
&&\hspace*{146pt}{} +
(1+\varepsilon^{-1}) \frac{1+\mathcal{N}^2/\ell
_0}{T^2}x^2 \biggr],\nonumber
\end{eqnarray}
inequality which holds on $\mathcal{B}$ with probability larger than
$1- 2\#\{
\mathcal{M}_{T}
\}e^{-x}$.
\item
Let us control now $A_2=2\nu_T(s_{\hat{m}}-s_m)$. To do so, we
need to control all the $V_{m'}=\nu_T(s_{m'}-s_m)$.
But on $\mathcal{B}$, $ V_{m'} = U_{m'}$ where
\[
U_{m'} = \frac{1}{T} \int_0^T \mathbh{1}_{t\leq\tau} \Psi_{s_{m'}-s_m}(t)
\bigl(dN_t-\Psi_{s}(t)\,dt\bigr).
\]
So one can use Corollary 1 of \cite{reynaudspl}: with probability
larger than $1-e^{-x}$,
\[
U_{m'} \leq\sqrt{2 v x} + \frac{b}{3}x,
\]
where $v$ and $b$ are constants such that for all $t\leq T$,
\[
v \geq\frac{1}{T^2}\int_0^T \mathbh{1}_{t\leq\tau} \Psi_{s_{m'}-s_m}(t)^2
\Psi_{s}(t) \,dt \quad\mbox{and}\quad b \geq\mathbh{1}_{t\leq\tau} \frac
{1}{T}\bigl|\Psi_{(s_{m'}-s_m)}(t)\bigr|.
\]

By similar arguments, we can obtain the following upper bound (see
\cite{sopacomplet}):
on $\mathcal{B}$ with probability larger than $1-\#\{\mathcal
{M}_{T}\} e^{-x}$
%
%
\begin{equation}
\label{majonu}
\nu_T(s_{\hat{m}}-s_m) \leq\Vert s_{\hat{m}}-s_m \Vert\sqrt{2
\frac
{(\eta
+H\mathcal{N})
R^2}{T} x} + \frac{2H\mathcal{N}}{3T} x.
\end{equation}
But $\Vert s_{\hat{m}}-s_m \Vert\leq
\Vert s_{\hat{m}}-s \Vert+\Vert s-s_m \Vert$. Thus, with the same constant
$\theta
$ as in (\ref{majoa}),
this gives (see \cite{sopacomplet})
%
%
\begin{equation}
\label{majob}
A_2 \leq\theta\Vert s_{\hat{m}}-s \Vert^2 + \theta\Vert s_m-s
\Vert^2 + \biggl(
\frac{4}{\theta}+\frac{2}{3R^2} \biggr)\frac{(\eta+ H\mathcal{N})
R^2}{T} x.
\end{equation}
\end{enumerate}
Now let us go back to (\ref{depart}). Using (\ref{majoa}) and (\ref
{majob}), we have actually obtained that
on $\mathcal{B}$ and on an event $\Omega_x$ whose probability is
larger than
$1-3\#\{\mathcal{M}_{T}\}e^{-x}$, the following inequality is true:
\begin{eqnarray*}
D^2_T(\tilde{s}-s) &\leq& D^2_T(s_m-s) +
\theta[\Vert\tilde{s}-s_{\hat{m}} \Vert^2+\Vert s_{\hat{m}}-s
\Vert^2]+
\theta
\Vert s-s_m \Vert^2 \\
&&{}+ \frac{1}{\theta}
\biggl[ (1+\varepsilon) (\eta+ H \mathcal{N})
R^2\frac{|\hat{m}|+1}{T} \bigl(1+3\sqrt{2x} \bigr)^2\\
&&\hspace*{95pt}{} +
(1+\varepsilon^{-1}) \frac{1+\mathcal{N}^2/\ell
_0}{T^2}x^2 \biggr]  \\
&&{} + \biggl( \frac{4}{\theta}+\frac{2}{3R^2} \biggr)\frac{(\eta+H\mathcal{N})
R^2}{T} x
+ \operatorname{pen}(m) - \operatorname{pen}(\hat{m}).
\end{eqnarray*}
As $s_\perp$ denotes the orthogonal projection for $\Vert\cdot \Vert$ of
$s$ on
$\mathcal{S}$, we can remark that
\[
\Vert\tilde{s}-s_{\hat{m}} \Vert^2+\Vert s_{\hat{m}}-s \Vert
^2=\Vert\tilde{s}-s \Vert^2=\Vert\tilde{s}-s_\perp\Vert
^2+\Vert s_\perp-s \Vert^2.
\]
Moreover,
\begin{eqnarray*}
D_T^2(\tilde{s}-s_\perp)&=&\frac1{T}\int_0^T \bigl(\Psi_{\tilde
{s}-s}(t)+\Psi_{s-s_\perp}(t)\bigr)^2 \,dt \\
&\leq&
(1+\varepsilon) D_T^2(\tilde
{s}-s) + (1+\varepsilon
^{-1})D_T^2(s-s_\perp).
\end{eqnarray*}
Hence, we obtain that on $\mathcal{B}\cap\Omega_x$
\begin{eqnarray*}
D^2_T(\tilde{s}-s_\perp) &\leq& (1+\varepsilon)
D^2_T(s_m-s) + (1+\varepsilon^{-1})
D_T^2(s-s_\perp)\\
&&{} + (1+\varepsilon)
\theta[\Vert\tilde{s}-s_\perp\Vert^2+\Vert s_\perp-s \Vert^2]
\\
&&{} +(1+\varepsilon) \theta
\Vert s-s_m \Vert^2 + (1+\varepsilon)\operatorname
{pen}(m)\\
&&{} + (1+\varepsilon) \biggl[\frac{1}{\theta}
(1+\varepsilon) (\eta+ H \mathcal{N})
R^2\frac{|\hat{m}|+1}{T} \bigl(1+3\sqrt{2x} \bigr)^2 -\operatorname{pen}(\hat
{m}) \biggr]\\
&&{} + (1+\varepsilon) \biggl( \frac{4}{\theta}+\frac{2}{3R^2}
\biggr)\frac{(\eta+H\mathcal
{N})R^2}{T} x\\
&&{} + \frac{(1+\varepsilon)
(1+\varepsilon^{-1})}{\theta} \frac{1+\mathcal
{N}^2/\ell_0}{T^2}x^2.
\end{eqnarray*}
But on $\mathcal{B}$, $D^2_T(\tilde{s}-s_\perp)\geq r^2\Vert\tilde
{s}-s_\perp\Vert^2$ since $\tilde{s}-s_\perp$ belongs to $\mathcal{S}$.
Hence, if we choose $\theta=r^2(1+\varepsilon)^{-2}$,
we obtain
\begin{eqnarray*}
\hspace*{-4pt}&&\frac{\varepsilon r^2}{1+\varepsilon}
\Vert\tilde{s}-s_\perp\Vert^2 \mathbh{1}_{\mathcal{B}\cap
\Omega_x} \\
\hspace*{-4pt}&&\qquad\leq(1+\varepsilon) D^2_T(s_m-s) +
(1+\varepsilon^{-1}) D_T^2(s-s_\perp)+ (1+\varepsilon)
\theta\Vert s_\perp-s \Vert^2\\
\hspace*{-4pt}&&\qquad\quad{} +(1+\varepsilon) \theta
\Vert s-s_m \Vert^2 + (1+\varepsilon)\operatorname
{pen}(m)\\
\hspace*{-4pt}&&\qquad\quad{} + (1+\varepsilon) \biggl( \frac{4}{\theta}+\frac{2}{3R^2}
\biggr)\frac{(\eta+H\mathcal
{N})R^2}{T} x+ \frac{(1+\varepsilon)
(1+\varepsilon^{-1})}{\theta} \frac{1+\mathcal
{N}^2/\ell_0}{T^2}x^2.
\end{eqnarray*}
It remains to add $\varepsilon r^2(1+\varepsilon)^{-1}\Vert s_\perp-s
\Vert^2 \mathbh{1}_{\mathcal
{B}\cap\Omega_x}$ on both sides, to obtain
(\ref{BESOIN}). For~(\ref{BESOIN2}),
let us take the expectation on both parts. We can remark that
$\mathbb{E}(D^2_T(s_m-s))=\Vert s_m-s \Vert_D^2\leq K^2 \Vert s_m-s
\Vert^2$,
by applying Lemma \ref{norme} and similar computations hold for
$s_\perp
$. Moreover, remark that
$\Vert s-s_\perp\Vert\leq\Vert s_m-s \Vert$,
since $S_m$ is a subset of~$\mathcal{S}$. This concludes the proof.
\end{pf*}
\begin{remark}\label{remark4}
In case of self-inhibition [see (\ref{sineg}) and Remark
\ref{remark3}], it is sufficient to replace $\mathcal{T}$ by $\mathcal
{T}\cap
\mathcal{T}'$ where
\[
\mathcal{T}'=\{t / \lambda(t)=0\}
\]
and to define accordingly the stopping time $\tau$ to obtain (\ref{BESOIN}).
\end{remark}

\subsection[Proof of Proposition 5]{Proof of Proposition \protect\ref{controldeb}}\label{sec83}
The control of $\mathcal{B}$ is twofold.

On one hand, one needs to control the number of points in any interval
of length~$A$. The control of the number of points in one interval comes
from some tedious computations that have been done in \cite{RB-Roy}.
Then the control for any interval comes from a reasoning that is close
in essence to the control of the suprema of identically distributed
variables with exponential moment.

On the other hand, one needs to control the deviations of $D_T^2(f)$
from its mean for $f$ in a finite vectorial subspace. We decompose the
problem in controlling the deviations of the associated bilinear form
for elements of the basis. Those deviations are controlled by using a
concentration inequality for Hawkes processes that have been derived via
coupling in \cite{RB-Roy}.

The heart of the proof actually consists in the probabilistic results
derived in~\cite{RB-Roy}. The final step is composed of lengthy and not
very informative computations that are omitted here and which can be
found in \cite{sopacomplet}.
\subsection[Proof of the minimax results (Propositions 2, 3 and 4)]{Proof of
the minimax results (Propositions \protect\ref{minimaxun}, \protect\ref{lowerholder} and
\protect\ref{minimaxislirr})}\label{sec84}

We first need two important lemmas.
\begin{Lemma}\label{kullback}
Let $f=(\mu,g)$ and
$s=(\nu,h)$ be two elements of $\mathbb{L}^{2}$ such that $\mu,\nu
>0$, $ g,h
\geq
0$, $\int g <1$ and $\int h <1$.
Let $\mathbb{P}_f^{[-A,T]}$, respectively, $\mathbb{P}_s^{[-A,T]}$, be
the distribution
of a stationary Hawkes process
with intensity $\Psi_{f}(\cdot)$, respectively, $\Psi_{s}(\cdot)$, restricted to
$[-A,T]$. Then the Kullback--Leibler distance satisfies
\begin{eqnarray*}
\mathbb{K}\bigl(\mathbb{P}_f^{[-A,T]},\mathbb{P}_s^{[-A,T]}\bigr)&=&\mathbb
{E}_f \biggl(\int_0^T
\phi\biggl[\log\biggl(\frac{\Psi_{s}(t)}{\Psi_{f}(t)} \biggr) \biggr]
\Psi_{f}(t) \,dt \biggr)\\
&&{} + \mathbb{K}\bigl(\mathbb{P}_f^{[-A,0]},\mathbb{P}_s^{[-A,0]}\bigr),
\end{eqnarray*}
where $\phi(u)=e^u-u-1$ and $\mathbb{E}_f$ represents the expectation with
respect to $\mathbb{P}_f^{[-A,T]}$.

Moreover, if $f$ and $s$ belong to $ \mathcal{L}_{H,P}^{\eta,\rho}$
and if
$A\Vert h \Vert_\infty\leq P-\log P-1$, then
\[
\mathbb{K}\bigl(\mathbb{P}_f^{[-A,T]},\mathbb{P}_s^{[-A,T]}\bigr)\leq T
\mathcal{C}_1 \Vert f-s \Vert^2 +
\mathcal{C}_2,
\]
where $\mathcal{C}_1$ and $\mathcal{C}_2$ are positive constants
depending only
on $A,H,P,\eta,\rho$.
\end{Lemma}

Lemma \ref{kullback} shows that the Kullback--Leibler distance between
two different processes linearly increases with $T$. It also clarifies
the link between the natural Kullback--Leibler distance and the
$\mathbb
{L}^2$-norm, $\Vert\cdot \Vert$, we used.
\begin{pf*}{Proof of Lemma \ref{kullback}}
Let us denote by $\mathbb{P}_f^{[0,T]}|_{[-A,0]}$ the conditional
distribution of the points of the process lying in $[0,T]$
conditionally to
the family of points lying in $[-A,0]$. Then the classical
decomposition of
the Kullback--Leibler distance with respect to the marginals gives the
following decomposition:
\[
\mathbb{K}\bigl(\mathbb{P}_f^{[-A,T]},\mathbb{P}_s^{[-A,T]}\bigr)=\mathbb
{E}_f \biggl[\ln\frac{d\mathbb{P}
_f^{[0,T]}|
_{[-A,0]}}{d\mathbb{P}_s^{[0,T]}|_{[-A,0]}} \biggr]+\mathbb{K}\bigl(\mathbb{P}
_f^{[-A,0]},\mathbb{P}
_s^{[-A,0]}\bigr).
\]
Next, we combine Example 7.2(b) with Proposition 7.2.III of
\cite{daleyVerejones} to obtain that the conditional likelihood ratio is
\[
\frac{d\mathbb{P}_f^{[0,T]}|_{[-A,0]}}{d\mathbb{P}_s^{[0,T]}|
_{[-A,0]}}=\exp
\biggl(\int_0^T
\ln[\Psi_{f}(t)/\Psi_{s}(t)] \,dN_t-\int_0^T \Psi_{f}(t)\,dt+\int_0^T
\Psi_{s}(t)\,dt \biggr).
\]
Using the martingale properties and the fact that the intensity is
predictable, one gets the first equation of Lemma \ref{kullback}.
Now to upper bound the Kullback--Leibler distance, we need first to
remark that
$\forall x>-1, \log(1+x)\geq x/(1+x)$ which gives that
\begin{eqnarray*}
\mathbb{E}_f \biggl(\int_0^T
\phi\biggl[\log\biggl(\frac{\Psi_{s}(t)}{\Psi_{f}(t)} \biggr) \biggr]
\Psi_{f}(t) \,dt \biggr)&\leq&
\mathbb{E}_f \biggl(\int_0^T\frac{(\Psi_{s}(t)-\Psi_{f}(t))^2}{\Psi
_{s}(t)}\,dt \biggr)\\
&\leq&
\frac
{T}{\rho} \Vert f-s \Vert_D^2.
\end{eqnarray*}
It is important to note that here (and only here) $\Vert\cdot \Vert_D$ is
computed with respect to $f$ and not $s$.
Now it remains to use Lemma \ref{norme} and to upperbound the constants
depending on $f$ by constants depending on $A,H,P,\eta,\rho$ to
obtain the
first part of the inequality.

Then it remains to upper bound $\mathbb{K}(\mathbb
{P}_f^{[-A,0]},\mathbb{P}_s^{[-A,0]})$.
This quantity is just a remaining term: we only need to prove that on
$\mathcal{L}^{\eta,\rho}_{H,P}$, this term cannot explode. A~lengthy
but necessary proof of it can be found in \cite{sopacomplet}. In
essence, it is close to Proposition \ref{controldeb} and it heavily
depends on the results of \cite{RB-Roy}.
\end{pf*}

Lemma \ref{kullback} combined with Birg\'e's lemma \cite{fanolucien}
gives the following result, which is ready to use for the different
lower bounds in the different situations.
\begin{Lemma}\label{fanopournous}
Let $\mathcal{S}$ be a family of possible $s$ such that $\Psi_{s}(\cdot)$ is
the intensity of a stationary Hawkes process, and such that $s$ belongs to
$\mathcal{L}_{H,p}^{\eta,\rho}$. Let $\delta>0$ and
let $\mathcal{C}\subset
\mathcal{S}$ be a finite family such that for all $f=(\mu,g)\in
\mathcal{C}$,
$A\Vert g \Vert_\infty\leq P-\log P-1$. Then there exists $\zeta_1$ and
$\zeta_2$ two particular positive functions of
$\eta,\rho,A,P,H$ such that if for all $f\not= f'$ in $\mathcal{C}$
\[
\frac{\zeta_1\log|\mathcal{C}|-\zeta_2}{T}\geq\Vert f-f' \Vert
^2\geq
\delta\qquad
\mbox{then }
\inf_{\hat{s}}\sup_{s\in\mathcal{S}} \mathbb{E}_s(\Vert\hat
{s}-s \Vert^2)\geq\frac
{\delta(1-\alpha)}{4},
\]
where $\alpha$ is an absolute positive constant (see \cite{fanolucien}
for a precise value).
\end{Lemma}
\begin{pf}
First, it is very classical to obtain that
\[
\inf_{\hat{s}}\sup_{s\in\mathcal{S}} \mathbb{E}_s(\Vert\hat
{s}-s \Vert^2)\geq
\frac14\inf_{\hat{s}\in\mathcal{C}}\sup_{s\in\mathcal{C}}
\mathbb{E}_s(\Vert\hat{s}-s \Vert^2).
\]
But
\[
\mathbb{E}_s(\Vert\hat{s}-s \Vert^2)\geq\delta\mathbb{P}_s(\hat
{s}\not= s).
\]
So
\[
\inf_{\hat{s}}\sup_{s\in\mathcal{S}} \mathbb{E}_s(\Vert\hat
{s}-s \Vert^2)\geq
\frac\delta4\inf_{\hat{s}\in
\mathcal{C}} \Bigl(1-\inf_{s\in\mathcal{C}}\mathbb{P}_s(\hat{s}= s) \Bigr).
\]
It remains to apply Birg\'e's lemma \cite{fanolucien}, by
upper bounding the mean Kullback--Leibler distance on
$\mathcal{C}$. Using Lemma \ref{kullback}, it remains only to choose
$\zeta_1$
and $\zeta_2$ according to $\mathcal{C}_1$ and $\mathcal{C}_2$.
This concludes the proof.
\end{pf}

It is now sufficient to apply the previous lemma for good choices of
$\mathcal{C}$.
\begin{pf*}{Proof of Proposition \ref{minimaxun}}
Let $m$ be a model. We set
$D=|m|$. Let
$\mathcal{P}_0$ be the maximal collection of subsets of $m$, such that
for all $\mathcal{I}\not=\mathcal{I}'$ in $\mathcal{P}_0$,
$|\mathcal{I}\Delta\mathcal{I}'|\geq\theta|m|$, then by \cite{gallager},
one has that $\log|\mathcal{P}_0|\geq\sigma|m|$, for $\theta$ and
$\sigma$ some
absolute constants.

Let
\[
\mathcal{C}_0= \biggl\{f_\mathcal{I}= \biggl(\rho, \sum_{I\in\mathcal
{I}}\frac{\varepsilon
}{\sqrt{\ell(I)}} \mathbh{1}_I \biggr),
\mathcal{I}\in\mathcal{P}_0 \biggr\},
\]
where $\varepsilon$ is a positive real number that will
be chosen later. To
ensure that $\mathcal{C}_0\subset\mathcal{L}_{H,P}^{\eta,\rho}$,
we need
that $\varepsilon\leq\min(H, P/A)\sqrt{\ell_0}$.
Moreover, to apply Lemma
\ref{fanopournous}, we need that $\varepsilon\leq
(P-\log P-1) \sqrt{\ell_0}/A$.

Now, for all $f_\mathcal{I},f_{\mathcal{I}'}$ in $\mathcal{C}_0$,
\[
\Vert f_\mathcal{I}-f_{\mathcal{I}'} \Vert^2=|\mathcal{I}\Delta
\mathcal
{I}'|\varepsilon^2
\geq\theta D \varepsilon^2.
\]

Moreover,
\[
\Vert f_\mathcal{I}-f_{\mathcal{I}'} \Vert^2 \leq\varepsilon^2 D.
\]

Finally, taking
\[
\varepsilon^2=\min\biggl(\frac{(\zeta_1D-\zeta_2)\sigma
}{TD}, \ell_0\min\bigl(H,
P/A,(P-\log P-1)/A\bigr)^2 \biggr),
\]
and applying Lemma \ref{fanopournous} gives the result.
\end{pf*}
\begin{pf*}{Proof of Proposition \ref{minimaxislirr}}
Let $\Gamma$ be a partition of $(0,A]$ and let us concentrate first on
the Islands set. Let $\mathcal{P}_1$ be the maximal collection of
subsets of $\Gamma$ with cardinal $D$, such that for all
$\mathcal{I}\not=\mathcal{I}'$ in $\mathcal{P}_1$,
$|\mathcal{I}\Delta\mathcal{I}'|\geq\theta D$, then by the Appendix
of~\cite{ptrfpois}, one has that $\log|\mathcal{P}_1|\geq\sigma D\log\frac
ND$, for $\theta$ and $\sigma$ some absolute constants. Let
\[
\mathcal{C}_1= \biggl\{f_\mathcal{I}= \biggl(\rho, \sum_{I\in\mathcal
{I}}\frac{\varepsilon
}{\sqrt{\ell(I)}} \mathbh{1}_I \biggr),
\mathcal{I}\in\mathcal{P}_1 \biggr\}.
\]
Then the same computations as before give the result for the Islands
set. But note that the set $\mathcal{C}_1$ is also included in
$S_{\Gamma,(2D+1)}^{\mathrm{irr}}$. Consequently,
the lower bound is also valid up to some multiplicative constant for
$S_{\Gamma,(2D+1)}^{\mathrm{irr}}$.
\end{pf*}
\begin{pf*}{Proof of Proposition \ref{lowerholder}}
For the H\"olderian family, let $\varphi$ be a positive continuous
function on $\mathbb{R}$, null outside $(0,A]$ and such that for all
$x,y\in\mathbb{R}$, $|\varphi(x)-\varphi(y)|\leq|x-y|^a$. Remark that
a quantity that only depends on $\varphi$ actually depends on $A$ and
$a$.

Let $m$ be a regular partition of $(0,A]$ in $D$ pieces. Let
$\varphi_D(x)= L D^{-a} \varphi(Dx)$. Let
$\mathcal{P}_0$ be defined as before and
\[
\mathcal{C}_2= \biggl\{s_{\mathcal{I}}=\biggl(\rho,\sum_{I\in\mathcal{I}}
\varphi_D(x-u_I)\biggr), \mathcal{I}\in\mathcal{P}_0 \biggr\},
\]
where $u_I$ is
the left extremity of $I$. To ensure that $\mathcal{C}_2\subset
\mathcal{L}_{H,p}^{\eta,\rho}$ and that $\Vert g \Vert_\infty\leq
(P-\log
P-1)/A$, we need that $D\geq c(A,a,H,P) L^{1/a}$, for some positive
continuous function $c$.

But for all $s_{\mathcal{I}},s_{\mathcal{I}'}$ in $\mathcal{C}_2$,
\[
\Vert s_\mathcal{I}-s_{\mathcal{I}'} \Vert^2=|\mathcal{I}\Delta
\mathcal{I}'|L^2
D^{-2a-1} \int\varphi^2 \geq\theta L^2
D^{-2a }\int\varphi^2.
\]

Moreover,
\[
\Vert s_\mathcal{I}-s_{\mathcal{I}'} \Vert^2\leq L^2
D^{-2a }\int\varphi^2.
\]

But note that for $D$ large enough $\zeta_1\sigma D-\zeta_2\geq\zeta'
D$ for
some other constant $\zeta'$.

It remains to choose
\[
D=\lozenge_{H,P,A,\rho,\eta,a}\max\bigl[(TL^2)^{1/(2a+1)},L^{1/a} \bigr]
\]
to obtain the result.
\end{pf*}

\section{Conclusion}\label{sec9}
We proposed a method based on model selection principle for
Hawkes processes that is proved
to be adaptive minimax with respect to certain classes of
functions.
In practice, the multiplicative constant in the penalty is calibrated
in a data-driven way that is proved to work well on simulations. In
particular, we designed a new method---namely the \textit{Islands
strategy} coupled with the angle penalty---that seems to be really
adapted to our biological
problem, namely characterizing the dependence between the occurrences
of a biological signal.
Moreover, it allows us to estimate the right range of
interaction.

This work asks, however, for several future developments. First, it
is
necessary to treat interaction with another type
of events (e.g., promoter/genes)
with the \textit{Islands strategy}. Next, a test procedure should be
applied to know whether the function $h$ is really nonzero.
This would be equivalent to testing whether there exists an interaction
or not.

\section*{Acknowledgments}
We would like to warmly thank Pascal Massart for his support, but also
Gaelle Gusto for a preliminary work during her Ph.D. thesis and Olivier
Catoni for his advice on the Kullback--Leibler distance. We also
thank the anonymous referees for their smart advice and careful
reading.

\printaddresses

\end{document}